\documentclass{amsart}
\usepackage{amsmath,amssymb, amsbsy}
\usepackage[dvips]{graphicx}
\usepackage{epsf}
\usepackage{color,graphics}
\newcommand{\media}{\mkern12mu\hbox{\vrule height4pt depth-3.2pt
    width6pt} \mkern-17.9mu\int} 
\newcommand{\mediapiccola}{\mkern5mu\hbox{\vrule height4pt depth-3.5pt
    width4pt} \mkern-12.5mu\int}\newcommand{\R}{{\mathbb R}} 
\newcommand{\N}{{\mathbb N}}
\newcommand{\Z}{{\mathbb Z}}
\newcommand{\weakly}{\rightharpoonup}
\renewcommand{\a }{\alpha }

\newcommand{\D }{\Delta }
\newcommand{\e }{\varepsilon }
\newcommand{\g }{\gamma}

\renewcommand{\l }{\lambda }
\renewcommand{\L }{\Lambda }
\newcommand{\n }{\nabla }

\newcommand{\Di}{{\mathcal D}^{1,2}(\R^N)}
\newcommand{\Dia}{{\mathcal D}_{a_{\lambda}}^{1,2}(\R^N)}
\newcommand{\dive }{\mathop{\rm div} }
\newcommand{\SO }{\mathop{\mathbb{SO}} }
\newcommand{\Dik}{{\mathcal D}^{1,2}_k(\R^N)}
\newcommand{\Dis}{{\mathcal D}^{1,2}_{\rm circ}(\R^N)}

\newenvironment{pf}{\noindent{\bf Proof}.\enspace}{\rule{2mm}{2mm}\medskip}
\newenvironment{pfn}[1]{\noindent{\bf Proof of {#1}.\enspace}}{\rule{2mm}{2mm}
\hfill\medskip}
\newtheorem{Theorem}{Theorem}[section]

\newtheorem{Lemma}[Theorem]{Lemma}
\newtheorem{Proposition}[Theorem]{Proposition}

\newtheorem{remark}[Theorem]{Remark}
\begin{document}

\title[Multi-polar symmetric Schr\"odinger equations]{Nonlinear Schr\"odinger
  equations with symmetric  multi-polar potentials}

\author[V. Felli \and S. Terracini]{Veronica Felli \and
  Susanna Terracini}

\address{Universit\`a di Milano Bicocca, Dipartimento di Matematica e
  Applicazioni, Via  Cozzi 53, 20125 Milano.}
\email{{\tt veronica.felli@unimib.it, susanna.terracini@unimib.it}. }

\thanks{Supported by Italy MIUR, national project
``Variational Methods and Nonlinear Differential Equations''.
}

\maketitle

\section{Introduction}
\noindent
Schr\"odinger equations with Hardy-type singular potentials have been
the object of a quite large interest in the recent literature, see e.g.
\cite{AFP, CCP, CTV, egnell, FG, GP, Jan,  RW,  SM, terracini}. The
singularity of inverse square potentials $V(x)\sim{\l}{|x|^{-2}}$ is
critical both from the mathematical and the physical  point of view. As it does not
belong to the Kato's class, it cannot be regarded as a lower
order perturbation of the laplacian but strongly influences the
properties of the associated 
Schr\"odinger operator. Moreover, from the point of view of nonrelativistic
quantum mechanics, among potentials of type
$V(x)\sim\l|x|^{-\a}$, the inverse square case represent a
{\em{transition threshold}}: for $\l<0$ and $\a>2$ (attractively singular
potential), the energy is not lower-bounded and a particle near the
origin in the presence  a potential of this type ``falls''
to the center, whereas if $\a<2$ the discrete spectrum has a
lower bound (see \cite{LL}). 
Moreover inverse square singular potentials
arise in many fields, such as quantum mechanics, nuclear physics, 
molecular physics, and quantum cosmology; we refer to \cite{FLS} 
for further discussion and motivation.

The case of multi-polar Hardy-type potentials was considered in \cite{FT,
  duyckaerts}. In particular in \cite{FT} the authors studied the
ground states of the following 
class of nonlinear elliptic 
equations with a critical power-nonlinearity and a potential
exhibiting multiple inverse square singularities:
\begin{equation}\label{eq:69}
\begin{cases}
-\D
v-{\displaystyle{\sum_{i=1}^k}}\dfrac{\l_i}{|x-a_i|^2}\,v=v^{2^*-1},\\[15pt]
v>0\quad\text{in }\R^N\setminus\{a_1,\dots,a_k\}.
\end{cases}
\end{equation}
For Schr\"odinger operators $-\Delta+V$ the potential term $V$
describes the interactions of the quantum particles with the
environment. Hence, multi-singular 
inverse-square potentials are associated with the interaction of particles
 with a finite number of electric dipoles. 
The mathematical interest in this problem rests in its
criticality, for the exponent of the nonlinearity as well as the
singularities share the same order of homogeneity with the laplacian.

The analysis 
carried out in \cite{FT} highlighted how the existence of 
solutions to (\ref{eq:69}) heavily depends on the strength  and the location 
of the singularities. For the scaling properties of the problem, the
mutual interaction among the poles actually depends only on the shape
of their configuration. 
When the poles form a symmetric structure, it is natural to wonder how
the symmetry affects such mutual interaction. 
The present paper means to study this aspect from the point of view of
 the existence of solutions inheriting the same symmetry properties as
 the set of singularities. 
More precisely we deal with a class of nonlinear elliptic
equations on $\R^N$, $N\geq 3$, involving a critical
power-nonlinearity  as well as a Hardy-type potential 
which is singular on sets exhibiting some simple kind of symmetry, as
depicted in figures 1--3.  

Let us start by considering a potential featuring 
multiple inverse square singularities located on the vertices of $k$-side 
regular  concentric polygons.
Let us write $\R^N=\R^2\times\R^{N-2}$. For $k\in\N$, we consider the
group $\Z_k\times\SO(N-2)$ acting on $\Di$ as 
$u(y',z)\to v(y',z)=u(e^{2\pi \sqrt{-1}/k}y',Tz)$, $T$ being any rotation of $\R^{N-2}$.

Let us consider $m$ regular polygons (with $k$ sides) centered at the
origin and lying on the plane $\R^2\times \{0\}\subset\R^N$. Let us denote as
$a_i^\ell$, $i=1,2,\dots,k$, the vertices of 
the $\ell$-th polygon, $\ell=1,2,\dots,m$. Since the polygons are
regular, we have $a_i^\ell=e^{2\pi \sqrt{-1}/k} a_{i-1}^\ell$.
We look for $\Z_k\times\SO(N-2)$-invariant  
solutions to the following equation 
\begin{equation}\label{eq:22}
\begin{cases}
-\D
v-\dfrac{\l_0}{|x|^2}v-{\displaystyle{\sum_{\ell=1}^m}}
{\displaystyle{\sum_{i=1}^k}}\dfrac{\l_{\ell}}
{|x-a_i^{\ell}|^2}\,v=v^{2^*-1},\\[15pt]
v>0\quad\text{in }\R^N\setminus\{0,a_i^{\ell}:\ 1\leq\ell\leq m,\ 
  1\leq i\leq k\},
\end{cases}
\end{equation}
where $2^*=\frac{2N}{N-2}$, i.e. for solutions belonging to  the space 
$$
\Dik=\{u(z,y)\in\Di:\ u(e^{2\pi \sqrt{-1}/k}z,y)=u(z,|y|)\},
$$
see figures 1--3.
Here $\Di$ denotes the closure  of the space 
${\mathcal D}(\R^N)$ of smooth functions with compact
support  with respect to the Dirichlet norm
$$
\|u\|_{\Di}:=\bigg(\int_{\R^N}|\n u|^2\,dx\bigg)^{1/2}.
$$
\noindent Let us denote as
$$
r_{\ell}=|a_1^{\ell}|=|a_2^{\ell}|=\dots=|a_k^{\ell}|, \quad\text{for
} \ell=1,2,\dots,m,
$$
the radius of the $\ell$-polygon and as
$$\Lambda_{\ell}=k\l_{\ell}, \quad \ell=1,\dots,m,$$
the total mass of poles located on the $\ell$-th polygon.
The
Rayleigh quotient associated with problem~(\ref{eq:22}) in $\Dik$ is
\begin{align}\label{eq:25}
S_k(\l_0,\L_1,\dots, \L_m)=\inf_{\substack{u\in\Dik\\
    u\not\equiv 0}} 
\frac{{\displaystyle{\int_{\R^N} |\n
      u|^2dx-
\int_{\R^N}\bigg(
\frac{\l_0}{|x|^2}+\sum_{\ell=1}^m\sum_{i=1}^k
  \frac{\L_{\ell}}{k|x-a_i^{\ell}|^2}\bigg)
u^2(x)\,dx}}}
{{\displaystyle{ 
    \bigg(\int_{\R^N}  |u|^{2^*}dx\bigg)^{2/2^*}}}}.
\end{align}
It is well known  that minimizers of (\ref{eq:25}) solve equation
(\ref{eq:22}) up to a Lagrange multiplier. Theorems \ref{t:ach} and
\ref{t:ach1} give sufficient conditions for the attainability of
$S_k(\l_0,\L_1,\dots, \L_m)$ for large values of $k$. 
Letting $k\to\infty$, the Schr\"odinger operator converges, in the
sense of distributions, to the operator associated with a continuous
distribution of  mass
  on concentric circles. We stress that the convergence
 of the potentials does not hold in the natural way, i.e. in $L^p_{\rm
   loc}(\R^N)$ for any $p\leq \frac N2$, because of the singularity. 
To formulate the limiting problem,
 for any $r>0$, we denote as 
$$
S_r:=\{(x,0)\in\R^2\times \R^{N-2}:\ |x|=r\}
$$
 the circle of radius $r$ lying on the plane $\R^2\times\{0\}\subset\R^N$
 and consider the distribution $\delta_{S_r}\in{\mathcal D}'(\R^N)$
 supported in $S_r$ and defined by
$$
\Big\langle_{\hskip-31pt{\mathcal D}'(\R^N)} \,\,\delta_{S_r}, \varphi
\Big\rangle_{{\mathcal D}(\R^N)} := \media 
_{S_r}\varphi(x)\,d\sigma(x)= \frac1{2\pi r}\int_{S_r}\varphi(x)\,d\sigma(x)
 \quad\text{for any }\varphi \in{\mathcal D}(\R^N),
 $$
where $d\sigma$ is the line element on $S_r$.  
We look for solutions to the following equation 
\begin{equation}\label{eq:55}
\begin{cases}
-\D
v-\dfrac{\l_0}{|x|^2}v-{\displaystyle{\sum_{\ell=1}^m}}
{\L_{\ell}}
\bigg(\delta_{S_{r_{\ell}}}*\dfrac1{|x|^2}\bigg)v=v^{2^*-1},
\\[15pt]
v>0\quad\text{in }\R^N\setminus\Big\{0, \bigcup_{ 1\leq \ell\leq
  m}S_{r_{\ell}}\Big\}, 
\end{cases}
\end{equation}
which are invariant by the action of the group $\SO(2)\times\SO(N-2)$. 
To this purpose, the natural space to set the problem is 
$$
\Dis=\{u(z,y)\in\Di:\ u(z,y)=u(|z|,|y|)\},
$$
and to consider the associated Rayleigh quotient
\begin{align}\label{eq:56}
S_{\rm circ}&(\l_0,\L_1,\dots,
\L_m)\\
\nonumber&=\inf_{\substack{u\in\Dis\\u\not\equiv 0}}
\frac{{\displaystyle{\int_{\R^N} |\n
      u|^2\,dx-\int_{\R^N}\bigg(\frac{\l_0}{|y|^2}+
\sum_{\ell=1}^m
  {\L_{\ell}}\media_{S_{r_{\ell}}}\frac{d\sigma(x)}
{|x-y|^2}\bigg)u^2(y)\,dy
}}}{{\displaystyle{ 
    \bigg(\int_{\R^N}  |u|^{2^*}dx\bigg)^{2/2^*}}}}.
\end{align}
The following theorem contains a 
Hardy type inequality for potentials which are singular at circles.

\begin{Theorem}\label{t:hardys-ineq}
Let $N\geq 3$ and $r>0$. For any $u\in\Di$ the map $y\mapsto
u(y)\int_{S_r}\frac{d\sigma(x)}{|x-y|^2}$ belongs to $L^2(\R^N)$ and 
\begin{equation}\label{eq:53}
\Big(\frac{N-2}2\Big)^2\int_{\R^N}|u(y)|^2\bigg(\media_{S_r}
\frac{d\sigma(x)}{|x-y|^2}\bigg)\,dy\leq\int_{\R^N}|\n u(y)|^2\, dy.
\end{equation}
Moreover the constant $\big(\frac{N-2}2\big)^2$ is optimal and not attained.
\end{Theorem}

\noindent Hardy type inequalities involving singularities at smooth
compact boundaryless manifolds have been considered by several authors, see
\cite{DD,FMT} and references therein. In the aforementioned papers,
the potentials taken into account are of the type $|\!\mathop{\rm
   dist}(x,\Sigma)|^{-2}$, where $ \mathop{\rm
   dist}(x,\Sigma)$ denotes the distance from a smooth compact
 manifold $\Sigma$. We point out that such kind of potentials are
 quite different from the ones we are considering. Indeed an explicit
 computation yields 
\begin{equation}\label{eq:64}
V^r(y):=\media_{S_r}
\frac{d\sigma(x)}{|x-y|^2}=\frac{1}{\sqrt{(r^2+|y|^2)^2-4r^2|y'|^2}}
\quad\text{for all } y=(y',z)\in \R^N=\R^2\times\R^{N-2}.
\end{equation}
Hence
\begin{align*}
V^r(y)\sim\frac1{|y|^2}\quad\text{as}
\quad|y|\to+\infty
\end{align*}
whereas
\begin{align*}
V^r(y)=\frac1{\sqrt{r^2+|y|^2-2r|y'|}}
\,\frac1{\sqrt{r^2+|y|^2+2r|y'|}}\sim\frac1{2r\big||y|-r\big|}
\quad\text{as}\quad  \mathop{\rm dist}(y,S_r)\to 0.
\end{align*}
Hence the singularity at the circle of $V^r$ is weaker 
than the inverse square distance potential considered in
\cite{DD,FMT}, but  has the same behavior at $\infty$. We also
remark that $V^r$ is ``regular'' in the sense of the
classification of singular potentials given in \cite{FLS}.

\medskip\noindent 
Arguing as in \cite[Proposition 1.1]{FT}, it is easy to verify that
solvability of equations (\ref{eq:22}) and (\ref{eq:55}) requires the
positivity of the associated quadratic forms. Let us consider for
example  the quadratic form
associated with potentials singular on circles, i.e.
$$
Q^{\rm circ}_{\l_0,\L_1,\dots,\L_m}=
\int_{\R^N} |\n
      u|^2\,dx-\int_{\R^N}\bigg(\frac{\l_0}{|y|^2}+
\sum_{\ell=1}^m
  {\L_{\ell}}\media_{S_{r_{\ell}}}\frac{d\sigma(x)}
{|x-y|^2}\bigg)u^2(y)\,dy.
$$
From (\ref{eq:53}) and Sobolev's inequality, it follows that
$$
Q^{\rm circ}_{\l_0,\L_1,\dots,\L_m}\geq
\bigg[1-\frac4{(N-2)^2}\big(\l_0^++{\textstyle{\sum}}_{\ell=1}^m\L_{\ell}^+\big)\bigg]
\int_{\R^N} |\n
      u|^2\,dx,
$$
where
$t^+:=\max\{t,0\}$ denotes the positive part.
Hence $Q^{\rm circ}_{\l_0,\L_1,\dots,\L_m}$ is positive definite
whenever 
\begin{equation}\label{eq:66}
\l_0^++\sum_{\ell=1}^m\L_{\ell}^+<\frac{(N-2)^2}4,
\end{equation}
see \cite[Proposition 1.2]{FT} for further discussion on the
positivity of the quadratic form. Condition~(\ref{eq:66}) also ensures
the positivity of the quadratic form associated to problem (\ref{eq:22}).

\medskip \noindent 
The attainability of $S_k(\l_0,\L_1,\dots, \L_m)$ and 
$S_{\rm circ}(\l_0,\L_1,\dots,\L_m)$ requires a delicate balance
between the contribution of positive and negative masses. In particular, if $N\geq 4$ and all
the masses have the same sign, $S_k(\l_0,\L_1,\dots, \L_m)$ and 
$S_{\rm circ}(\l_0,\L_1,\dots,\L_m)$ are never achieved.
\begin{Theorem}\label{t:posmas}
Let $N\geq4$,  $\l_0,\Lambda_1,\dots,\Lambda_m\in\R$,
$r_1,r_2,\dots,r_m\in\R^+$ satisfy (\ref{eq:66}). If 
\begin{align*}
\text{either}\quad{\rm(i)}&\quad \L_{\ell}<0\hskip1.655cm\text{for all }\ell=1,\dots,m\\
\text{or}\quad{\rm(ii)}&\quad \l_0>0,\
\L_{\ell}>0\quad\text{for all }\ell=1,\dots,m
\end{align*}
then neither $S_k(\l_0,\L_1,\dots, \L_m)$ nor 
$S_{\rm circ}(\l_0,\L_1,\dots,\L_m)$ are attained.
\end{Theorem}

\medskip \noindent 
The analysis we are going to carry out in the present paper  will
highlight that,  from the point of    
view of minimization of the Rayleigh quotient, spreading mass 
 all over a continuum is more convenient than localization of mass at
 isolated points. 
\begin{Theorem}\label{t:achcirc}
Let $N\geq 4$, $\l_0,\Lambda_1,\dots,\Lambda_m\in\R$,
$r_1,r_2,\dots,r_m\in\R^+$ satisfy (\ref{eq:66}) and
\begin{align}\label{eq:57}
&\sum_{\ell=1}^m\Lambda_{\ell}\leq 0,\qquad
\l_0<\frac{(N-2)^2}4,
\end{align}
and 
\begin{align}\label{eq:58}
\begin{cases}
{\displaystyle{\sum_{\ell=1}^m\frac{\Lambda_{\ell}}{|r_{\ell}|^2}>0}},
&\text{if}\quad 
\l_0\leq\dfrac{N(N-4)}4,\\[10pt]
{\displaystyle{\sum_{\ell=1}^m
\frac{\Lambda_{\ell}}
{|r_{\ell}|^{\sqrt{(N-2)^2-4\l_0}}}>0}}, &\text{if}\quad
\dfrac{N(N-4)}4<\l_0<\dfrac{(N-2)^2}4.
\end{cases}
\end{align}
Then  the infimum in (\ref{eq:56}) is
achieved. In particular equation 
(\ref{eq:55}) admits a solution which is $\SO(2)\times\SO(N-2)$-invariant.
\end{Theorem}

As problem (\ref{eq:56}) is the limit of (\ref{eq:25}),   
 when $k\to\infty$, we expect the assumptions of Theorem
 \ref{t:achcirc} to  ensure the existence of solutions to (\ref{eq:25}) provided $k$ is sufficiently large.
Indeed the theorem below states that \eqref{eq:57} and \eqref{eq:58} are sufficient conditions on radii and masses
of the polygons for the infimum in
(\ref{eq:25}) to be achieved when $k$ is large.

\begin{Theorem}\label{t:ach}
Let $N\geq 4$, $\l_0,\Lambda_1,\dots,\Lambda_m\in\R$,
$r_1,r_2,\dots,r_m\in\R^+$ satisfy (\ref{eq:66}), \eqref{eq:57} and \eqref{eq:58}. 
For any $\ell=1,\dots,m$ and $k\in\N$ let $\{a_i^\ell\}_{i=1,\dots,k}$ be the vertices of a
regular $k$-side polygon centered at $0$ of radius $r_{\ell}$ and let
$\l_{\ell}=\Lambda_{\ell}/k$. Then if $k$ is sufficiently large,
 the infimum in (\ref{eq:25}) is achieved. In particular equation
(\ref{eq:22}) admits a solution which is $\Z_k\times\SO(N-2)$-invariant.
\end{Theorem}

When $N>4$, it is possible to estimate how large $k$ must be in order to obtain the above existence result. This is the content of the following theorem.

\begin{Theorem}\label{t:ach1}
Assume that $N>4$,  $\l_0^++k\sum_{j=1}^{m}\l_j^+<\frac{(N-2)^2}4$,
\begin{align}
\label{eq:42}
&\sum_{\ell=1}^m\l_{\ell}\leq 0,\\
\label{eq:43}
&\l_1\leq\l_2\leq\dots\leq\l_m\leq\frac{N(N-4)}4,\\
\label{eq:45}
&\l_0<\frac{(N-2)^2}4,\\
\label{eq:47}
&\begin{cases}
{\displaystyle{\sum_{\ell=1}^m\frac{\l_{\ell}}{|r_{\ell}|^2}>0}},
&\text{if}\quad 
\l_0\leq\dfrac{N(N-4)}4,\\
{\displaystyle{\sum_{\ell=1}^m
\frac{\l_{\ell}}
{|r_{\ell}|^{\sqrt{(N-2)^2-4\l_0}}}>0}}, &\text{if}\quad
\dfrac{N(N-4)}4<\l_0<\dfrac{(N-2)^2}4,
\end{cases}\\
\label{eq:46}
&\frac{\l_0}{|r_m|^2}+\l_m
\sum_{i=1}^{k-1}\frac{1}
{4r_m^2\big|\sin\frac{i\pi}k\big|^2}+\sum_{\ell=1}^{m-1}\l_{\ell}\sum_{i=1}^k\frac{1}{r_m^2+r_{\ell}^2-2 r_m
r_{\ell}\cos\Big(\frac{2\pi i}k+\Theta_{j\ell}\Big)
}>0,
\end{align}
where $\Theta_{j\ell}$ denoted the minimum angle formed by vectors
$a_i^j$ and $a_s^{\ell}$, (see figure 4).
Then the infimum in (\ref{eq:25}) is achieved. In particular equation
(\ref{eq:22}) admits a solution.
\end{Theorem}

This paper is organized as follows. In section
\ref{sec:problem-with-one} we recall some known facts about the
single-polar problem
and study the behavior of any solution (radial and non radial) to the
one pole-equation  
near the singularities $0$ and $\infty$. In section
\ref{sec:hardys-ineq-with} we prove  
the Hardy type inequality for potentials which are singular at circles
stated in Theorem \ref{t:hardys-ineq}. Section
\ref{sec:pala-smale-cond} contains an analysis of possible reasons for lack
of compactness of minimizing sequences of problems (\ref{eq:25}) and
(\ref{eq:56}) and a local Palais-Smale
condition below some critical thresholds.  In section
\ref{sec:inter-estim} we provide some interaction estimates which are
needed in section \ref{sec:comp-betw-conc} to compare the
concentration levels of minimization sequences and consequently to
prove Theorem \ref{t:achcirc}. Section \ref{sec:limit-s_kl-as}
contains the study of behavior of energy levels of minimizing
sequences as $k\to\infty$ which is needed in section
\ref{sec:proof-theor-reft} to prove Theorem \ref{t:ach}. Last section
is devoted to the proof of Theorem \ref{t:ach1}. 

\medskip
Figures 1--3 are the plot of the test functions used to estimate the
energy levels of Palais-Smale sequences and represent the expected
shape of  solutions  found in
Theorem \ref{t:ach} (assuming the
last $N-2$ variables to be $0$), which is  based on the knowledge of
their behavior at singularities (see 
\cite{FS3} and Theorem \ref{t:caratt}), which is known to be be
singular at poles with positive mass and vanishing at 
poles with negative mass.

\begin{center}
 \vskip0.5truecm\noindent
 \epsfxsize=5.5in \epsfbox{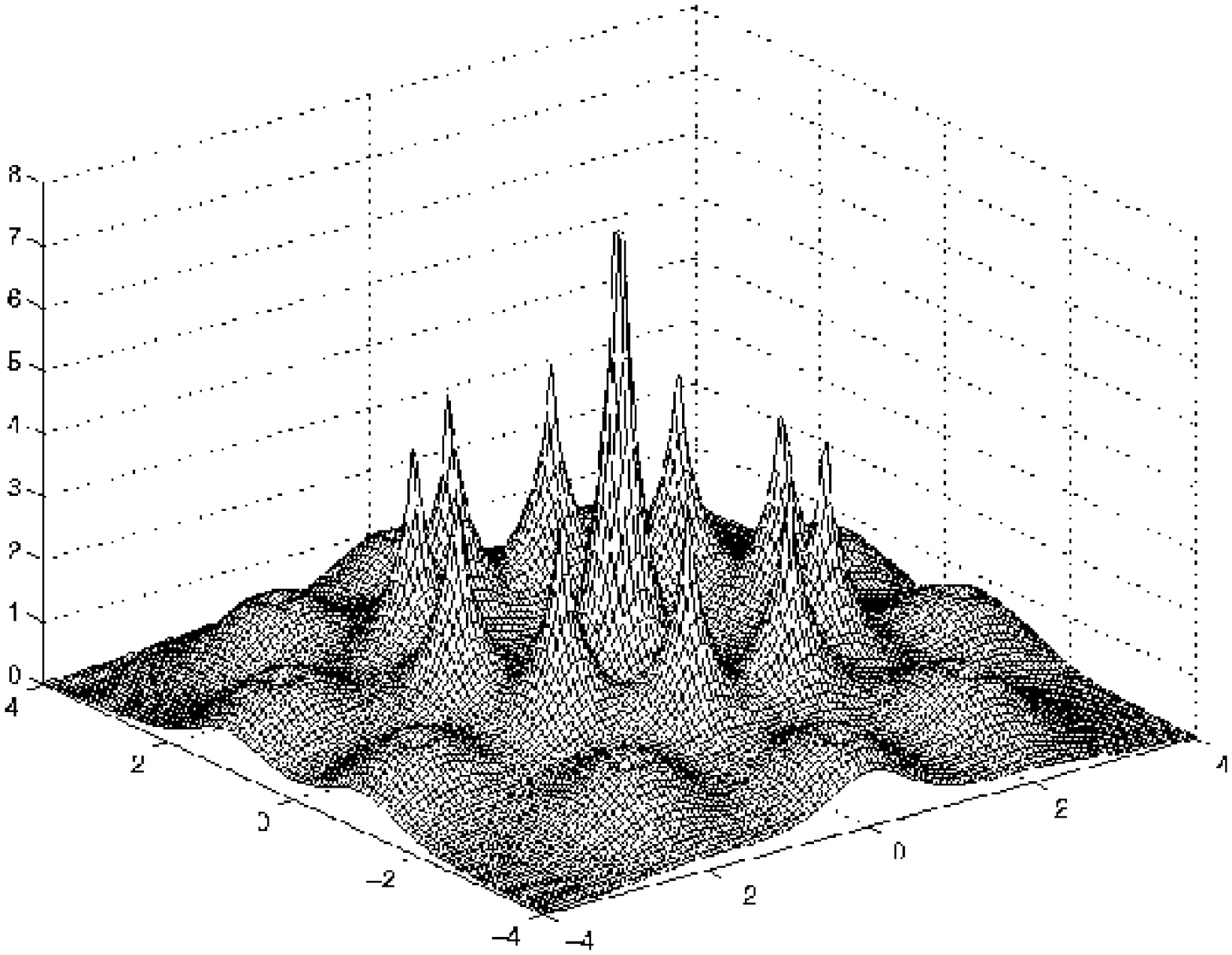}
\\
 {\scriptsize \ \ {\bf  1.}  $k=10$, $m=2$, $\l_0,\l_1>0$, $\l_2<0$.}
\vskip0.5truecm\noindent
\end{center}

\noindent
{\bf Notation. } We list below  some notation used throughout the
paper.\par
\begin{itemize}
\item[-]$B(a,r)$ denotes the ball $\{x\in\R^N: |x|<r\}$ in $\R^N$ with
  center at $a$ and radius $r$.
\item[-] $S_r=\{(x,0)\in\R^2\times \R^{N-2}: |x|=r\}$ denotes the circle
  of radius $r$ in the plane $\R^2\times\{0\}$.
\item[-] $\delta_x$ denotes the Dirac mass located at point $x\in\R^N$.
\item[-] ${\mathcal D}(\R^N)$ is the space of smooth functions with compact
support in $\R^N$.  
\item[-] $\Di$ is the closure of ${\mathcal D}(\R^N)$ with respect to the Dirichlet norm
$(\int_{\R^N}|\n u|^2\,dx)^{1/2}$.
\item[-] $C_0(\R^N)$ denotes the
closure of continuous functions with compact support in $\R^N$ with
respect to the uniform norm.
\item[-] $\mathop{\rm dist}(x,A)$ denotes the distance of the  point
  $x\in\R^N$ from the set $A\subset\R^N$.
\item[-] $\|\cdot\|_p$ denotes the norm in the Lebesgue space $L^p(\R^N)$.
\item[-] $A\triangle B$ denotes the symmetric difference of sets $A$
  and $B$, i.e. $A\triangle B=(A\setminus B)\cup (B\setminus A)$.

\begin{center}
 \vskip0.5truecm\noindent
 \begin{tabular}{c}
 \leavevmode
 \epsfxsize=5in \epsfbox{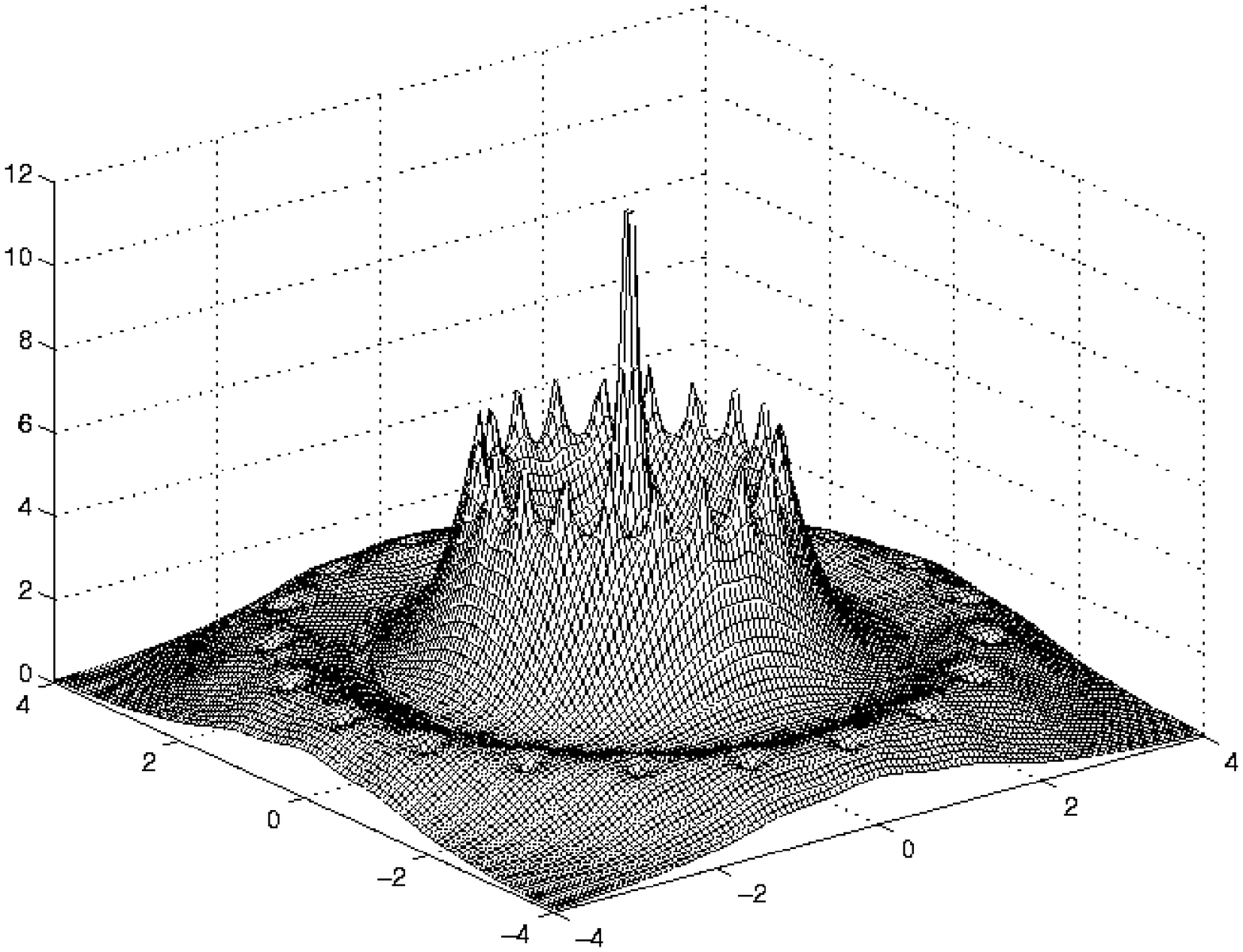} \\
\\
 {\scriptsize \ \ {\bf  2.} $k=20$, $m=2$, $\l_0,\l_1>0$, $\l_2<0$.}\\
 \epsfxsize=5in \epsfbox{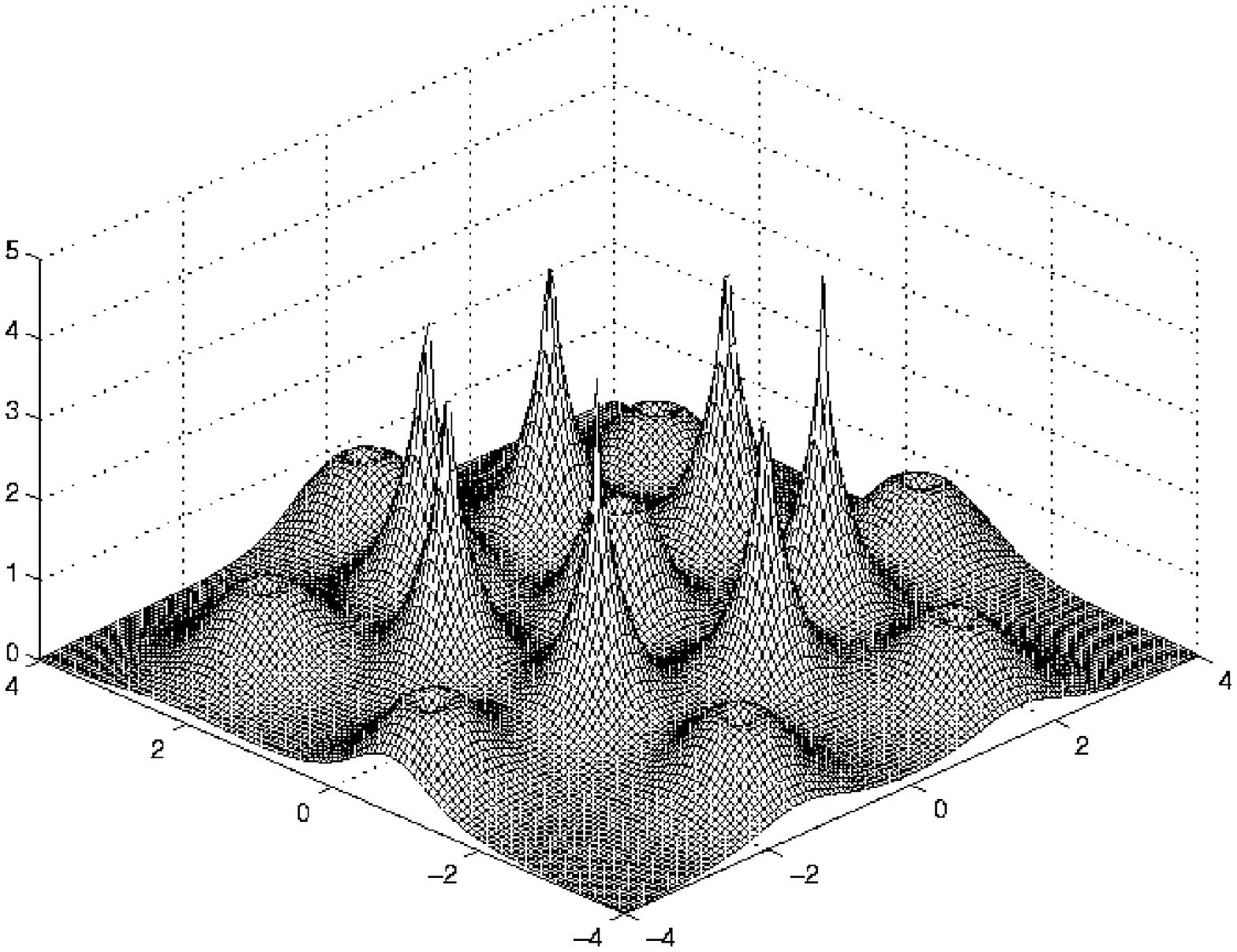} \\ 
 \\
{\scriptsize \ \ {\bf  3.} $k=6$, $m=2$, $ \l_1>0$, $\l_0,\l_2<0$.}
  \end{tabular}\\
\vskip0.5truecm\noindent
\end{center}

\end{itemize}

\section{The problem  with one singularity}\label{sec:problem-with-one}

\noindent For any $\l<(N-2)^2/4$, the  problem with one singularity  
\begin{equation}\label{eq:1}
\begin{cases}
-\D u=\dfrac{\l}{|x|^2}u+u^{2^*-1},\quad x\in \R^N,\\[10pt]
 u>0 \hbox{  in  }\R^N\setminus\{0\}, \text{  and  }u\in \Di,
\end{cases}
\end{equation}
 admits a family of  positive solutions
given by 
\begin{equation}\label{eq:ZA}
 w_{\mu }^{\l}(x)=\mu
^{-\frac{N-2}{2}}w^{\l}_1\Big(\dfrac{x}{\mu }\Big),\,\,
\mu>0,
\end{equation}
where we denote
\begin{equation}\label{eq:zzz}
w^{\l}_1(x)=\frac{\left(N(N-2)\nu
_\l^2\right)^{\frac{N-2}{4}}}{\left(|x|^{1-\nu _\l}(1+|x|^{2\nu
_\l})\right)^\frac{N-2}{2}},\quad \mbox{and}\quad\nu
_\l=\left(1-\frac{4\l}{(N-2)^2}\right)^{1/2}.
\end{equation}
Moreover, when $0\leq \l<(N-2)^2/4$,  all $w_{\mu }^{\l}(x)$ minimize the associated
Rayleigh quotient and the minimum can be computed as:
\begin{equation}\label{eq:minimiz}
S(\l):=\inf_{u\in \Di\setminus
\{0\}}\dfrac{Q_\l(u)}{\big(\int_{\R^N}  |u|^{2^*}dx\big)^{2/2^*}}
=\dfrac{Q_\l(w^{\l}_{\mu})}{\big(\int_{\R^N}  |w^{\l}_{\mu}|^{2^*}dx\big)^{2/2^*}}
=\bigg(1-\frac{4\l}{(N-2)^2}\bigg)^{\!\!\frac{N-1}{N}}S,
\end{equation}
where  we denoted the quadratic form $Q_\l(u)=\int_{\R^N} |\n u|^2dx -\l\int_{\R^N}
\frac{u^2}{|x|^2}dx$,
see
\cite{terracini}, and $S$ is  the best constant in the Sobolev
inequality
 \begin{equation*}
S\|u\|_{L^{2^*}(\R^N)}^2\leq \|u\|_{{\mathcal D}^{1,2}(\R^N)}^2.
\end{equation*}

\noindent As minimizers of problem (\ref{eq:minimiz}), we consider
\begin{equation}\label{eq:44}
z^{\l}_{\mu}(x)=\frac{w^{\l}_{\mu}(x)}{\big(\int_{\R^N}
  |w^{\l}_{\mu}|^{2^*}dx\big)^{1/2^*}}=\a_{\l,N}\,\mu^{-\frac{N-2}2}
\bigg(\Big|\frac{x}{\mu}\Big|^{1-\nu_{\l}}+
\Big|\frac{x}{\mu}\Big|^{1+\nu_{\l}} \bigg)^{-\frac{N-2}2}
\end{equation}
where $\a_{\l,N}=\left(N(N-2)\nu
_\l^2\right)^{\frac{N-2}{4}}\|w^{\l}_1\|_{L^{2^*}}^{-1}$
 is a positive constant depending only on $\l$ and $N$, so that
for $0\leq\l<(N-2)^2/4$
\begin{equation}\label{eq:62}
S(\l)=Q_{\l}(z^{\l}_{\mu})\quad\text{for all }\mu>0.
\end{equation}
For $-\infty<\l<(N-2)^2/4$, we also set
\begin{equation}\label{eq:minimizk}
S_k(\l):=\inf_{u\in \Dik\setminus
\{0\}}\dfrac{Q_\l(u)}{\big(\int_{\R^N}  |u|^{2^*}dx\big)^{2/2^*}}.
\end{equation}
We note that $S(\l)\leq S_k(\l)$, and equality holds whenever $\l\geq0$.
Moreover the following result has been proved in
\cite{terracini}.
\begin{Lemma}[see \cite{terracini}, Lemma 6.1]\label{l:ska}
Let $N\geq 4$. If $S_k(\l)<k^{2/N}S$ then $S_k(\l)$ is achieved.
\end{Lemma}

\medskip\noindent
 In \cite{terracini} it is proved that if $\l\in
(0,(N-2)^2/4)$ then all solutions to (\ref{eq:1}) are of the
form~(\ref{eq:ZA}) while if $\l<\!\!<0$ then also non radial solutions
to 
(\ref{eq:1}) can exist.
The behavior of any solution (radial and non radial) to
  problem (\ref{eq:1})
near the singularities $0$ and $\infty$ is described by the following 
theorem.
\begin{Theorem}\label{t:caratt}
If $\l<(N-2)^2/4$ and $u\in\Di$ is a solution to problem (\ref{eq:1}),
then there exist positive constant $\kappa_0(u)$ and $\kappa_{\infty}(u)$ 
 depending on $u$ such that
\begin{align}
\label{eq:3}
&u(x)=|x|^{-\frac{N-2}2(1-\nu_{\l})}\big[\kappa_0(u)+O(|x|^{\alpha})\big],\quad
\text{as }x\to 0,\\
\label{eq:4}
&u(x)=|x|^{-\frac{N-2}2(1+\nu_{\l})}
\big[\kappa_{\infty}(u)+O(|x|^{-\alpha})\big],\quad\text{as }|x|\to +\infty,
\end{align}
for some $\a\in(0,1)$.
\end{Theorem}

\begin{remark} Putting together (\ref{eq:3}--\ref{eq:4}) we deduce
  that there exists 
a positive constant $\kappa(u)$  depending on $u$ such that
\begin{equation}\label{eq:2}
\frac{1}{\kappa(u)}\,w^{\l}_1(x)\leq u(x)\leq \kappa(u)\,w^{\l}_1(x).
\end{equation}
\end{remark}

\begin{pfn}{Theorem \ref{t:caratt}}
Set 
$$
a_{\l}=\frac{(N-2)\big(1-\nu_{\l}\big)}2, \quad v(x)=|x|^{a_{\l}}u(x).
$$
Then the function $v$ belongs to $\Dia$ where
 $\Dia$ denotes the space obtained by completion of 
${\mathcal D}(\R^N)$ with respect to the weighted Dirichlet norm
$$
\|v\|_{\Dia}:=\bigg(\int_{\R^N}|x|^{-2 a_{\l}}|\n v|^2\,dx\bigg)^{1/2}.
$$
Moreover $v$ solves equation 
\begin{equation}\label{eq:6}
-\dive\big(|x|^{-2a_{\l}}\nabla v\big)=\frac{v^{2^*-1}}{|x|^{2^*
    a_{\l}}}.
\end{equation}
From \cite[Theorem 1.2]{FS3}, it follows that $v$ is H\"older
continuous; in particular expansion (\ref{eq:3}) holds for
$\kappa_0(u)=v(0)$ and some
$\alpha\in(0,1)$. Moreover $v(0)$ is strictly positive in view of 
Harnack's inequality for degenerate operators proved in \cite{Gutierrez}, see
also \cite{DV}; we mention that weights of type $|x|^{-2 a}$
with $a<\frac{N-2}2$ belong to the class of quasi-conformal weights
considered in \cite{Gutierrez}.

To deduce (\ref{eq:4}), we perform the change of
variable
$$
\tilde v(x)=|x|^{a_{\l}-(N-2)}u\bigg(\frac{x}{|x|^2}\bigg),
$$
and observe that the transformed function $\tilde v$ solves
equation (\ref{eq:6}). Hence \cite[Theorem 1.2]{FS3} yields that
$\tilde v$
is H\"older continuous and admits the following expansion
$$
\tilde v(x)=\tilde v(0)+O(|x|^{\alpha}),\quad
\text{as }x\to 0,
$$
for some $\alpha\in(0,1)$, where $\tilde v(0)>0$ in  view of 
Harnack's inequality  in \cite{Gutierrez}. Coming back to $u$ we
obtain that $u$ satisfies (\ref{eq:4}) with $\kappa_{\infty}(u)=\tilde
v(0)>0$. 
\end{pfn}

\section{Hardy's inequality with singularity on a
  circle}\label{sec:hardys-ineq-with} 
\noindent  We prove now the  
Hardy type inequality for potentials which are singular at circles
stated in Theorem \ref{t:hardys-ineq}.

\medskip
\begin{pfn}{Theorem \ref{t:hardys-ineq}}
Let us consider the minimization problem 
$$
I(S_r):=\inf_{\substack{u\in\Di\\u\not\equiv0}}\frac{{\displaystyle{\int_{\R^N}|\n u(y)|^2\, dy}}}
{{\displaystyle{\int_{\R^N}|u(y)|^2\bigg(\media_{S_r}
\frac{d\sigma(x)}{|x-y|^2}\bigg)\,dy}}}
=\inf_{\substack{u\in {\mathcal D}(\R^N\setminus\{0\})
\\u\not\equiv 0}}\frac{{\displaystyle{\int_{\R^N}|\n u(y)|^2\, dy}}}
{{\displaystyle{\int_{\R^N}|u(y)|^2\bigg(\media_{S_r}
\frac{d\sigma(x)}{|x-y|^2}\bigg)\,dy}}}
,
$$
where the last equality is due to density of ${\mathcal D}(\R^N\setminus\{0\})$ in 
$\Di$ (see e.g. \cite[Lemma 2.1]{catrinawang}). 
An easy calculation shows that for any $u\in\Di$
$$
\frac{{\displaystyle{\int_{\R^N}|\n u(y)|^2\, dy}}}
{{\displaystyle{\int_{\R^N}|u(y)|^2\bigg(\media_{S_r}
\frac{d\sigma(x)}{|x-y|^2}\bigg)\,dy}}}
=\frac{{\displaystyle{\int_{\R^N}|\n v(y)|^2\, dy}}}
{{\displaystyle{\int_{\R^N}|v(y)|^2\bigg(\media_{S_1}
\frac{d\sigma(x)}{|x-y|^2}\bigg)\,dy}}}
$$
where $v(y)=u(r y)$. Therefore 
\begin{equation}\label{eq:54}
I(S_r)=I(S_1)\quad\text{for any }r>0.
\end{equation}
In view of (\ref{eq:54}), it is enough to prove the theorem for $r=1$.
The proof consists in three steps.

\medskip\noindent
{\bf Step 1: Inequality (\ref{eq:53}) i.e. $\boldsymbol{I(S_1)\geq
    \big(\frac{N-2}2\big)^2}$.} 

\smallskip
For any $u\in \Di$, $u\geq 0$ a.e., we consider the Schwarz symmetrization $u^*$ of $u$ defined as 
\begin{equation}\label{eq:82}
u^*(x):=\inf\big\{t>0:\ \big|\{y\in\R^N:\ u(y)>t\}\big|\leq \omega_N|x|^N\big\}
\end{equation}
where $|\cdot|$ denotes the Lebesgue measure of $\R^N$ and $\omega_N$ is the volume of the standard unit $N$-ball. From  \cite[Theorem 21.8]{Willem}, it follows that for any $x\in S_1$
$$
\int_{\R^N}\frac{|u(y)|^2}{|x-y|^2}\,dy\leq \int_{\R^N}|u^*(y)|^2\left[\left(\frac1{|x-y|}\right)^*\right]^2.
$$
Since $\big(\frac1{|x-y|}\big)^*=\frac1{|y|}$, we deduce
\begin{equation}\label{eq:wil}
\int_{\R^N}\frac{|u(y)|^2}{|x-y|^2}\,dy\leq \int_{\R^N}\frac{|u^*(y)|^2}{|y|^2}\,dy.
\end{equation}
Moreover by Polya-Szego inequality
\begin{equation}\label{eq:wil1}
\int_{\R^N}|\n u^*|^2\leq 
\int_{\R^N}|\n u|^2.
\end{equation}
From (\ref{eq:wil}--\ref{eq:wil1}) and the classical Hardy's
inequality,  it follows that, for any $u\in
\Di\setminus\{0\}$, $u\geq 0$ a.e., 
\begin{align*}
\frac{{\displaystyle{\int_{\R^N}|\n u(y)|^2\, dy}}}
{{\displaystyle{\int_{\R^N}|u(y)|^2\bigg(\media_{S_1}
\frac{d\sigma(x)}{|x-y|^2}\bigg)\,dy}}}&=
\frac{{\displaystyle{\int_{\R^N}|\n  u(y)|^2\, dy}}}
{{\displaystyle{\media_{S_1}
\bigg(\int_{\R^N}| u(y)|^2\frac{dy}{|x-y|^2}\bigg)\,d\sigma(x)}}}\\[10pt]
&
\geq 
\frac{{\displaystyle{\int_{\R^N}|\n u^*(y)|^2\, dy}}}
{{\displaystyle{\int_{\R^N}\frac{|u^*(y)|^2}{|y|^2}\,dy}}}
\geq \Big(\frac{N-2}2\Big)^2.
\end{align*} 
Due to evenness of the quotient we are minimizing, to compute $I(S_1)$
it is enough to take the infimum over positive functions. Hence
passing to the infimum in the above inequality, we obtain 
$I(S_1)\geq \big(\frac{N-2}2\big)^2$.

\goodbreak
\medskip\noindent
{\bf Step 2: Optimality of the constant,
  i.e. $\boldsymbol{I(S_1)=\big(\frac{N-2}2\big)^2}$.}

\smallskip
We  fix $u\in {\mathcal D}(\R^N\setminus\{0\})$ and let $0<r<R$ be
such that $\mathop{\rm supp}u\subset\{x\in \R^N:\ r<|x|<R\}$.  
For any $0<\l<\!\!\!<1$, we set $\tilde u_{\l}(x)=u(\l x)$. Hence we have
\begin{align}\label{eq:step2}
I(S_1)&\leq  
\frac{{\displaystyle{\int_{\R^N}|\n \tilde u_{\l}(y)|^2\, dy}}}
{{\displaystyle{\int_{\R^N}|\tilde u_{\l}(y)|^2\bigg(\media_{S_1}
\frac{d\sigma(x)}{|x-y|^2}\bigg)\,dy}}}=
\frac{{\displaystyle{\int_{\R^N}|\n u(y)|^2\, dy}}}
{{\displaystyle{\int_{\R^N}|u(y)|^2\bigg(\media_{S_1}
\frac{d\sigma(x)}{|\l x-y|^2}\bigg)\,dy}}}.
\end{align}
Since $\frac{1}{|\l x-y|^2}\leq \frac4{r^2}$ for all $y\in \mathop{\rm supp}u$, $x\in S_1$, and $0<\lambda<\frac r2$, by Dominated Convergence Theorem we deduce that $\int_{\R^N}|u(y)|^2\big(\mediapiccola_{S_1}
\frac{d\sigma(x)}{|\l x-y|^2}\big)\,dy$ converges to 
 $\int_{\R^N}\frac{|u(y)|^2}{|y|^2}\,dy$  as $\l\to 0$, hence passing to the limit in \eqref{eq:step2} we obtain 
$$
I(S_1)\leq 
\frac{{{\int_{\R^N}|\n u(y)|^2\, dy}}}
{{{\int_{\R^N}|y|^{-2}|u(y)|^2\,dy}}}
\quad\text{for any\ } u\in u\in {\mathcal D}(\R^N\setminus\{0\}).
$$
By density of ${\mathcal D}(\R^N\setminus\{0\})$ in $\Di$ we deduce
$$
I(S_1)\leq 
\inf_{\Di\setminus\{0\}}
\frac{{
{\int_{\R^N}|\n u(y)|^2\, dy}}}
{{
{\int_{\R^N}|y|^{-2}|u(y)|^2\,dy}}}
=\Big(\frac{N-2}2\Big)^2
$$
where the last equality follows from the optimality of the constant
$\big(\frac{N-2}2\big)^2$ in the classical Hardy inequality (see
\cite[Lemma 2.1]{GP}. Collecting the above inequality with the one
proved in Step 1, we find $I(S_1)=\big(\frac{N-2}2\big)^2$.

\medskip\noindent
{\bf Step 3: The infimum $\boldsymbol{I(S_1)}$ is not attained.}

\smallskip

Arguing by contradiction, assume that the infimum $I(S_1)$ is achieved
by some $\bar u\in\Di$. We can assume that $\bar u\geq 0$ (otherwise
we consider $|\bar u|$ which is also a minimizer by evenness of the
quotient). Hence from \eqref{eq:wil} and \eqref{eq:wil1}
$$
\Big(\frac{N-2}2\Big)^2=
\frac{{\displaystyle{\int_{\R^N}|\n \bar u(y)|^2\, dy}}}
{{\displaystyle{\media_{S_1}
\bigg(\int_{\R^N}|\bar u(y)|^2\frac{dy}{|x-y|^2}\bigg)\,d\sigma(x)}}}
\geq 
\frac{{\displaystyle{\int_{\R^N}|\n \bar u^*(y)|^2\, dy}}}
{{\displaystyle{\int_{\R^N}\frac{|\bar u^*(y)|^2}{|y|^2}\,dy}}}
\geq \Big(\frac{N-2}2\Big)^2.
$$
Therefore the above inequalities are indeed equalities and this
implies that the infimum 
\begin{equation}\label{eq:chhar}
\Big(\frac{N-2}2\Big)^2=\inf_{u\in\Di\setminus\{0\}}
\frac{{\displaystyle{\int_{\R^N}|\n  u(y)|^2\, dy}}}
{{\displaystyle{\int_{\R^N}{|y|^{-2}}{| u(y)|^2}\,dy}}}
\end{equation}
which yields the best constant in the classical Hardy inequality, is
achieved by $\bar u^*$. Since it is known that the infimum in
\eqref{eq:chhar} cannot be attained (see \cite[Remark
1.2]{terracini}), we reach a contradiction.~\end{pfn}

\section{The  Palais-Smale condition under $\Z_k\times\SO(N-2)$ and $\SO(2)\times\SO(N-2)$-invariance}\label{sec:pala-smale-cond}

\noindent Let us define the functional $J_k:\Di\to\R$ associated to
equation (\ref{eq:22}) as 
\begin{align}\label{eq:energy} 
J_k(u)=&
\frac12\int_{\R^N} |\n u|^2dx-
\frac{\l_0}2\int_{\R^N}
\frac{u^2(x)}{|x|^2}\,dx\\
\nonumber&-
\sum_{\ell=1}^m\sum_{i=1}^k
\frac{\l_{\ell}}2\int_{\R^N}
\frac{u^2(x)}{|x-a_i^{\ell}|^2}\,dx-
\frac{S_k(\l_0,\L_1,\dots, \L_m)}{2^*}\int_{\R^N} |u|^{2^*}dx.
\end{align}
The choice of location of the singularities ensures that $J_k$ is 
$\Z_k\times\SO(N-2)$-invariant. Since $\Z_k\times\SO(N-2)$ acts by
isometries on $\Di$, we can apply the {\sl Principle of Symmetric
    Criticality} by Palais \cite{palais} to deduce that the critical
    points of $J_k$ restricted to $\Dik$ are also critical points of $J_k$
    in $\Di$.
Therefore, if $u$ is a critical point of $J_k$ in $\Dik$, $u>0$
outside singularities, 
then  $v=S_k(\l_0,\L_1,\dots,
\L_{\ell})^{1/(2^*-2)}u$ is a solution
to equation~(\ref{eq:22}).

The following theorem provides a  local Palais-Smale
condition for $J_k$ restricted to $\Dik$  below some critical threshold.
We emphasize that the invariance of the problem by the action of a 
subgroup of orthogonal transformation allows to recover some
compactness, in the sense that concentration points of invariant
functions  
must be located in some symmetric way, thus reducing the possibility
of loss of compactness. The restriction on dimension $N\geq 4$ is
required to avoid the presence of possible concentration points on
$\{0\}\times\R^{N-2}$. Indeed  when $N=3$, $\SO(N-2)=\SO(1)$ is a
discrete group, making thus possible concentration at
 points on the axis $\{0\}\times\R$.

We mention that the {\it
  Concentration-Compactness} 
method under the action of $\Z_k\times\SO(N-2)$ was used by several authors
 to find $k$-bump
solutions with prescribed symmetry for different classes of nonlinear
elliptic equations: nonlinear Schr\"odinger equation in \cite{wang},
nonlinear elliptic equations in symmetric domains in \cite{catrinawang3},
nonlinear elliptic equations of Caffarelli-Kohn-Nirenberg type in
\cite{catrinawang2}, and
elliptic equations with Hardy potential and critical growth in
\cite{terracini}.

\begin{Theorem}\label{th:ps}
Assume $N\geq 4$ and $\l_0^++k\sum_{j=1}^{m}\l_j^+<\frac{(N-2)^2}4$. Let $\{u_n\}_{n\in\N}\subset \Dik$
be a Palais-Smale sequence for
$J_k$ restricted to $\Dik$, namely 
$$ 
\lim_{n\to\infty}J_k(u_n)=c<\infty\text{ in }\R\quad\text{and}\quad
\lim_{n\to\infty}J_k'(u_n)=0\text{ in the dual space }(\Dik)^{\star}. 
$$ 
If 
\begin{align}\label{eq:19}
\!c<\frac{S_k(\l_0,\L_1,\dots,\L_{m})^{1-\frac{N}{2}}}
{N}
\min\!\bigg\{\!k^{\frac2N}S,
k^{\frac2N}S(\l_1), \dots, k^{\frac2N}S(\l_m), S_k(\l_0),
S_k\Big(\l_0+k\sum\limits_{\ell=1}^{m}\l_{\ell}\Big)\!\bigg\}^{\!\!\frac
  N2}\!\!\!, 
\end{align}
then $\{u_n\}_{n\in\N}$ has a subsequence strongly converging in $\Dik$.
\end{Theorem}
\begin{pf}
Let $\{u_n\}_{n\in\N}$ be a Palais-Smale sequence for $J_k$ in $\Dik$, then from
Hardy's and Sobolev's inequalities it is easy to prove that $\{u_n\}_{n\in\N}$
is a bounded sequence in $\Di$.  Hence, up to
a subsequence, $u_n\rightharpoonup u_0\mbox{ in } \Di$ and 
$ u_n\to u_0$  almost everywhere.
Therefore, from the {\sl Concentration Compactness Principle} by
P. L. Lions (see 
\cite{PL1, PL2}), we deduce the existence of a subsequence,
still denoted 
by $\{u_n\}$, for which there exist an at most
countable set ${\mathcal J}$, points $x_j\in\R^N\setminus\{0,a_i^{\ell},\
1\leq i\leq k,\ 1\leq \ell\leq m\}$,
real numbers $\mu_{x_j},\nu_{x_j}$, $j\in{\mathcal
  J}$, and $\mu_0$, $\nu_0$, $\mu_{a_i^{\ell}}$, $\nu_{a_i^{\ell}}$, $\g_{a_i^{\ell}}$,
$i=1,\dots,k$, $\ell=1,\dots,m$, such that the
following convergences hold in the sense 
of measures
\begin{align}
\label{eq:21}&|\nabla u_n|^2\rightharpoonup d\mu\ge |\n
u_0|^2+\mu_0\delta_0+\sum_{\ell=1}^m\sum_{i=1}^k \mu_{a_i^{\ell}}
\delta_{a_i^{\ell}}+\sum_{j\in {\mathcal
    J} }\mu_{x_j}\delta_{x_j},\\
\label{eq:27}&|u_n|^{2^*}\rightharpoonup d\nu=
|u_0|^{2^*}+\nu_0\delta_0+
\sum_{\ell=1}^m\sum_{i=1}^k
\nu_{a_i^{\ell}}\delta_{a_i^{\ell}}+ \sum_{j\in {\mathcal
    J}}\nu_{x_j}\delta_{x_j},\\
\label{eq:33}& \l_0\dfrac{u_n^2}{|x|^2}\rightharpoonup 
d \g_{0}=
\l_{0}\dfrac{u_0^2}{|x|^2}+\g_{0}
\delta_{0},\\
\label{eq:28}& \l_\ell\dfrac{u_n^2}{|x-a_i^{\ell}|^2}\rightharpoonup 
d \g_{a_i^{\ell}}=
\l_{\ell}\dfrac{u_0^2}{|x-a_i^{\ell}|^2}+\g_{a_i^{\ell}}
\delta_{a_i^{\ell}},\quad\text{for any}\quad
i=1,\dots,k,\quad \ell=1,\dots,m.
\end{align}
From Sobolev's inequality it follows that 
\begin{equation}\label{eq:30}
S\nu_{x_j}^{\frac2{2^*}}\le\mu_{x_j} \mbox{ for all } j\in{\mathcal
  J}\quad\text{and}\quad S\nu_{a_i}^{\frac2{2^*}}\le\mu_{a_i} \mbox{ for all }
i=1,\dots,k.
\end{equation}
Possible concentration at infinity of  the sequence can be quantified
by the following numbers 
\begin{equation}\label{eq:59}
\nu_{\infty}=\lim_{R\to \infty}\limsup_{n\to
\infty}\int_{|x|>R}|u_n|^{2^*}dx,\quad \mu_{\infty}=\lim_{R\to
\infty}\limsup_{n\to \infty}\int_{|x|>R}|\n u_n|^{2}dx,
\end{equation}
and  $$
\g_{\infty}=\lim_{R\to \infty}\limsup_{n\to
\infty}\int_{|x|>R}\Big(\l_0+k\sum_{\ell=1}^m\l_{\ell}\Big)\frac{ u_n^2}{|x|^2}dx.
$$
Let us first prove that possible concentration points are located in a
symmetric fashion.
Pointwise convergence of $u_n\in\Dik$ to $u_0$ implies that $u_0$ is
invariant by the $\Z_k\times\SO(N-2)$-action, hence $u_0\in\Dik$.
Moreover, for any $\phi\in C_0(\R^N)$ and for any $\tau\in\Z_k\times\SO(N-2)$
we have
$$
\int_{\R^N}|u_n|^{2^*}\phi=\int_{\R^N}|u_n|^{2^*}(\phi\circ\tau^{-1})
$$
that, passing to the limit as $n\to\infty$, yields
$$
\sum_{\ell=1}^m\sum_{i=1}^k\nu_{a_i^{\ell}}\phi(a_i^{\ell})+
\sum_{j\in{\mathcal J}}\nu_{x_j}\phi(x_j)=
\sum_{\ell=1}^m\sum_{i=1}^k\nu_{a_i^{\ell}}\phi(\tau^{-1}(a_i^{\ell}))+
\sum_{j\in{\mathcal J}}\nu_{x_j}\phi(\tau^{-1}(x_j)).
$$
Choosing $\phi=\phi_{i,\e}^{\ell}$ such that
$\phi_{i,\e}^{\ell}\equiv1$ on $B(a_i^{\ell},\e/2)$,
$\phi_{i,\e}^{\ell}\equiv0$ on $\R^N\setminus B(a_i^{\ell},\e)$,
$0\leq \phi_{i,\e}^{\ell}\leq 1$, and letting $\e\to 0$, we find that
$\nu_{a_i^{\ell}}= \nu_{\tau^{-1}(a_i^{\ell})}$. Hence we deduce that
for any $\ell=1,\dots,m$ there exists $\nu_a^{\ell}\in\R$ such that 
$$
\nu_{a_i^{\ell}}=\nu_a^{\ell}
\quad\text{for any }i=1,\dots,k.
$$
Fix $j\in{\mathcal J}$ and let  $\phi=\phi_{j,\e}$ such that
$\phi_{j,\e}\equiv1$ on $B(x_j,\e/2)$,
$\phi_{j,\e}\equiv0$ on $\R^N\setminus B(x_j,\e)$,
$0\leq \phi_{j,\e}\leq 1$. Letting  $\e\to 0$, we find that
\begin{align*}
{\rm(i)}&\qquad\text{either }\nu_{x_j}=0\\
{\rm(ii)}&\qquad\text{or for any }\tau\in Z_k\times\SO(N-2) \text{ there exists }
i\in{\mathcal J} \text{ such that } \tau(x_j)=x_i.
\end{align*}
When $N\geq 4$, $\SO(N-2)$ is a continuous group, hence for any
$x\not\in \R^2\times \{0\}\subset\R^N$, the set $\{\tau(x):\ \tau\in
Z_k\times\SO(N-2)\}$ is more than countable. 
If alternative (ii) holds, from at most countability of ${\mathcal J}$
we deduce that $x_j\in \R^2\times \{0\}\subset\R^N$.
Moreover, arguing as above we can prove that if $\tau(x_j)=x_i$ for
some $\tau\in Z_k\times\SO(N-2)$ then $\nu_{x_i}=\nu_{x_j}$.
Hence we can rewrite (\ref{eq:27}) as 
$$
|u_n|^{2^*}\rightharpoonup d\nu=
|u_0|^{2^*}+\nu_0\delta_0+
\sum_{\ell=1}^m \nu_a^{\ell}
\sum_{i=1}^k
\delta_{a_i^{\ell}}+ \sum_{\ell\in {\mathcal
    L}}\nu_{y}^{\ell}\sum_{i=1}^k\delta_{y_i^{\ell}}
$$
where ${\mathcal L}$ is an at most countable set,
$$
\{y_i^{\ell}:\ i=1,\dots,k,\ \ell\in {\mathcal L}\}\subset \{x_j:\
j\in{\mathcal J}\}\quad\text{and}\quad y_i^{\ell}=e^{2\pi
   \sqrt{-1}/k}y_{i-1}^{\ell}. 
$$
{\bf Concentration at non singular points.}\ We claim that
\begin{equation}\label{eq:5}
\mathcal J\quad\text{is finite and for }j\in {\mathcal J}\text{ 
either }\nu_{x_j}=0\text{ or }\nu_{x_j}\geq
\Big(\frac{S}{S_k(\l_0,\L_1,\dots, \L_m)}\Big)^{N/2}.
\end{equation}
Indeed, if $\nu_{x_j}>0$, we have $x_j=y^{\ell}_{s}$ and
$\nu_{x_j}=\nu_y^{\ell}$ for some $1\leq s\leq k$, $\ell\in{\mathcal L}$.
 For  $\e >0$, let $\phi_{\ell}^{\e}$ be a smooth  cut-off
function in $\Dik$, $0\leq \phi^{\e}_{\ell}(x)\leq 1$ such that
$$
 \phi_{\ell}^{\e}(x)= 1\quad\mbox{ if }x\in\bigcup_{i=1}^kB\Big(y_i^{\ell},
 \frac{\e }{2}\Big),\quad 
\phi_{\ell}^{\e}(x)= 0\quad\mbox{ if
}x\not\in\bigcup_{i=1}^kB\big(y_i^{\ell},\e \big),\quad\text{and}\quad
|\nabla 
\phi_{\ell}^{\e} |\le \dfrac{4}{\e }. 
$$
Testing $J_k'(u_n)$ with $u_n\phi_{\ell}^{\e}$ we obtain
\begin{align*}
 0&=\lim_{\e\to 0}\lim_{n\to \infty }\langle
J_k'(u_n),u_n\phi_{\ell}^{\e}\rangle \geq\sum_{i=1}^k \mu_{y_i^{\ell}}
-kS_k(\l_0,\L_1,\dots,\L_m)\nu_{y}^{\ell}.
\end{align*} 
By (\ref{eq:30}) we
have that $ S(\nu_{y}^{\ell})^{\frac2{2^*}}\le\mu_{y_i^{\ell}}$, then
we obtain that
$\nu_{y}^{\ell}\geq \big(\frac{S}{S_k(\l_0,\L_1,\dots,
\L_m)}\big)^{N/2}$. Hence 
either $\nu_{x_j}=0$ or $\nu_{x_j}\geq \big(\frac{S}{S_k(\l_0,\L_1,\dots,
\L_m)}\big)^{N/2}$, which implies that
${\mathcal J}$ is finite. In particular also ${\mathcal L}$ is finite.
The claim is proved. 

\medskip\noindent{\bf Concentration at vertices of polygons.}\ We claim that
\begin{equation}\label{eq:9}
\text{ for each}\quad \ell=1,2,\dots,m
\quad\text{either}\quad\nu_{a}^{\ell}=0\quad\text{or}\quad
\nu_{a}^{\ell}\geq
\Big(\frac{S(\l_{\ell})}{S_k(\l_0,\L_1,\dots,\L_m)}\Big)^{N/2}.
\end{equation}
In order to prove claim (\ref{eq:9}), for each $i=1,2,\dots,k$ and
$\ell=1,2,\dots,m$ we consider the smooth 
cut-off function  $\phi_{i,\e}^{\ell}$ satisfying $0\leq \phi_{i,\e}^{\ell}(x)\leq 1$,
$$
 \phi_{i,\e}^{\ell}(x)= 1\quad\mbox{ if }|x-a_i^{\ell}|\leq
 \frac{\e}{2},\quad  
\phi_{i,\e}^{\ell}(x)= 0\quad\mbox{ if }|x-a_i^{\ell}|\geq
\e,\quad\text{and}\quad |\nabla 
\phi_{i,\e}^{\ell} |\le \dfrac{4}{\e }. 
$$
From \eqref{eq:minimiz} we obtain
that
\begin{equation*}
\frac{\int_{\R^N} |\n (u_n\phi_{i,\e}^{\ell})|^2dx -\l_{\ell}\int_{\R^N}
{|x-a_i^{\ell}|^{-2}}{|\phi_{i,\e}^{\ell}|^2u^2_n}\,dx}{\Big(\int_{\R^N}|\phi_{i,\e}^{\ell}
u_n|^{2^*}\Big)^{2/2^*}}\geq
S(\l_{\ell})
\end{equation*}
hence passing to limit as $n\to\infty$ and $\e\to 0$ we obtain
\begin{equation}\label{eq:10}
\mu_{a_i^{\ell}}\geq \g_{a_i^{\ell}}+
S(\l_{\ell})\big(\nu_{a}^{\ell}\big)^{2/2^*}.
\end{equation}
For  $\e >0$, let $\psi_{\ell}^{\e}$ be a smooth  cut-off
function in $\Dik$, $0\leq \psi^{\e}_{\ell}(x)\leq 1$ such that
$$
 \psi_{\ell}^{\e}(x)= 1\quad\mbox{ if }x\in\bigcup_{i=1}^kB\Big(a_i^{\ell},
 \frac{\e }{2}\Big),\quad 
\psi_{\ell}^{\e}(x)= 0\quad\mbox{ if
}x\not\in\bigcup_{i=1}^kB\big(a_i^{\ell},\e \big),\quad\text{and}\quad
|\nabla 
\psi_{\ell}^{\e} |\le \dfrac{4}{\e }. 
$$
Testing $J_k'(u_n)$ with $u_n\psi_{\ell}^{\e}$ and letting $n\to\infty$ and $\e\to 0$ 
we infer that
\begin{equation}\label{eq:11}
\sum_{i=1}^k\mu_{a_i^{\ell}}-\sum_{i=1}^k\g_{a_i^{\ell}}\leq k S_k(\l_0,\L_1,\dots,\L_m)\nu_{a}^{\ell}. 
\end{equation}
From (\ref{eq:10}) and (\ref{eq:11}) we 
derive \eqref{eq:9}.

\medskip\noindent{\bf Concentration at the origin.}\ We claim that
\begin{equation}\label{eq:20}
\text{either}\quad\nu_0=0\quad\text{or}\quad
\nu_0\geq
\Big(\frac{S_k(\l_0)}{S_k(\l_0,\L_1,\dots,\L_m)}\Big)^{N/2}.
\end{equation}
In order to prove claim (\ref{eq:20}), we consider a smooth 
cut-off function  $\psi_{0}^{\e}\in \Dik$ satisfying $0\leq \psi_{0}^{\e}(x)\leq 1$,
$$
 \psi_{0}^{\e}(x)= 1\quad\mbox{ if }|x|\leq
 \frac{\e}{2},\quad  
\psi_{0}^{\e}(x)= 0\quad\mbox{ if }|x|\geq
\e,\quad\text{and}\quad |\nabla 
\psi_{0}^{\e} |\le \dfrac{4}{\e }. 
$$
From \eqref{eq:minimizk} we obtain
that
\begin{equation*}
\frac{\int_{\R^N} |\n (u_n\psi_{0}^{\e})|^2dx -\l_{0}\int_{\R^N}
{|x|^{-2}}{|\psi_{0}^{\e}|^2u^2_n}\,dx}{\Big(\int_{\R^N}|\psi_{0}^{\e}
u_n|^{2^*}\Big)^{2/2^*}}\geq
S_k(\l_0)
\end{equation*}
hence passing to limit as $n\to\infty$ and $\e\to 0$ we obtain
\begin{equation}\label{eq:34}
\mu_0\geq \g_0+
S_k(\l_0)\big(\nu_0\big)^{2/2^*}.
\end{equation}
On the other hand, testing $J_k'(u_n)$ with $u_n\psi_{0}^{\e}$ and
letting $n\to\infty$ and $\e\to 0$  
we infer that
\begin{equation}\label{eq:35}
\mu_{0}-\g_{0}\leq  S_k(\l_0,\L_1,\dots,\L_m)\nu_{0}. 
\end{equation}
From (\ref{eq:34}) and (\ref{eq:35}) we deduce (\ref{eq:20}).

\medskip\noindent{\bf Concentration at infinity.}\ We claim that
\begin{equation}\label{eq:12}
\text{either}\quad\nu_{\infty}=0\quad\text{or}\quad \nu_{\infty}\geq
\bigg(\frac{S_k(\l_0+k\sum_{\ell=1}^m
  \l_{\ell})}{S_k(\l_0,\L_1,\L_2,\dots,\L_m)}\bigg)^{{N}/{2}}. 
\end{equation}
In order to prove (\ref{eq:12}),
we  study the possibility of concentration at  $\infty$.
Let $\psi_R$ be a regular radial cut-off function such that 
$$
0\le \psi_R(x)\le 1,\quad
\psi_R(x)=
\left\{
\begin{array}{l}
1,\,\mbox{ if }|x|>2R \\ 0,\,\mbox{ if }|x|<R,
\end{array}
\right. \quad\text{and}\quad |\nabla \psi_R |\le \dfrac{2}{R} . 
$$
From \eqref{eq:minimizk} we obtain
that
\begin{equation}\label{eq:13}
\frac{\int_{\R^N} |\n (u_n\psi_R)|^2dx
  -\Big(\l_0+k\sum_{\ell=1}^m\l_{\ell}\Big)\int_{\R^N}
\frac{\psi_R^2u^2_n}{|x|^2}dx}{\Big(\int_{\R^N}|\psi_R
u_n|^{2^*}\Big)^{2/2^*}}\geq
S_k\Big({\textstyle{\l_0+k\sum_{\ell=1}^m\l_{\ell}}}\Big).
\end{equation}
Taking $\limsup$ as $n\to\infty$ and limit as $R\to+\infty$, standard
calculations yield
\begin{equation}\label{eq:15}
\mu_{\infty}-\g_{\infty}\geq
S_k\Big({\textstyle{\l_0+k\sum_{\ell=1}^m\l_{\ell}}}\Big)
\nu_{\infty}^{2/2^*}.
\end{equation} 
Testing $J_k'(u_n)$ with $u_n\psi_R$ and letting $n\to\infty$ and $R\to+\infty$, we obtain
\begin{equation}\label{eq:17}
\mu_{\infty}-\g_{\infty}\leq S_k(\l_0,\L_1,\dots,\L_m)\nu_{\infty}. 
\end{equation}
Claim (\ref{eq:12}) follows from (\ref{eq:15}) and (\ref{eq:17}).

\medskip\noindent
As a conclusion we obtain
\begin{align}\label{eq:18} 
c & =  J_k(u_n)-\frac{1}{2}\langle J_k'(u_n),u_n\rangle
  +o(1)=  \frac{1}{N}S_k(\l_0,\L_1,\dots,\L_m)\int_{\R^N}
|u_n|^{2^*}dx +o(1)\\
& = \frac{S_k(\l_0,\L_1,\dots,\L_m)}{N}\bigg\{\int_{\R^N}
|u_0|^{2^*}dx+\nu_0+k\sum_{\ell=1}^m\nu_{a}^{\ell}+k\sum_{\ell\in{\mathcal
    L}}\nu_y^{\ell}
+\nu_{\infty}\bigg\}.\notag
\end{align}
From (\ref{eq:19}), (\ref{eq:5}), (\ref{eq:9}), (\ref{eq:20}), (\ref{eq:12}), 
and (\ref{eq:18}), we deduce that $\nu_0=0$, $\nu_{y}^{\ell}=0$ for any $\ell\in {\mathcal
  L}$, $\nu_{a}^{\ell}=0$ for any $\ell=1,\dots,m$, and $\nu_{\infty}=0$. Then,
up to a subsequence, $u_n\to u_0$ in $\Dik$.
\end{pf}

\noindent The functional $J:\Di\to\R$ associated to
equation (\ref{eq:55}) is 
\begin{align}\label{eq:energy-circ} 
J(u)=&
\frac12\int_{\R^N} |\n u|^2dx-
\frac{\l_0}2\int_{\R^N}
\frac{u^2(x)}{|x|^2}\,dx\\
\nonumber&-
\sum_{\ell=1}^m
\frac{\L_{\ell}}2\int_{\R^N}\bigg(\media_{S_{r_{\ell}}}\frac{u^2(y)}
{|x-y|^2}\,d\sigma(x)\bigg)\,dy-
\frac{S_{\rm circ}(\l_0,\L_1,\dots, \L_m)}{2^*}\int_{\R^N} |u|^{2^*}dx.
\end{align}
The functional  $J$ is 
$\SO(2)\times\SO(N-2)$-invariant. Since $\SO(2)\times\SO(N-2)$ acts by
isometries on $\Di$, the {\sl Principle of Symmetric
    Criticality} by Palais \cite{palais} implies that  the critical
    points of $J$ restricted to $\Dis$ are also critical points of $J$
    in $\Di$. Therefore, any critical point of $J$ in $\Dis$ provides 
 a solution to equation~(\ref{eq:55}).

The following theorem is the analogous of Theorem \ref{th:ps} for  $J$
restricted to $\Dis$. However, the fact that the singularities are
spread over circles instead of being concentrated at atoms reduces the
possibility of lack of compactness. Indeed, according to P.L. Lions  {\sl
  Concentration-Compactness Principle}, possible forms of ``non
compactness'' which can cause failure of Palais-Smale condition are
loss of mass at infinity and concentration at an most countable set of
 points. When considering $J$ restricted to $\Dis$, it turns out that
 if $\bar x$ is a concentration point of a Palais-Smale sequence, then all
 points of the orbit ${\mathcal O}(\bar x)=\{\tau\,\bar x:\ \tau\in\SO(2)\times\SO(N-2)\}$
 must be concentration points. On the other hand, 
when $N\geq 4$ both groups $\SO(2)$ and $\SO(N-2)$ are continuous,
hence the only point $\bar  x$ for which ${\mathcal O}(\bar x)$ is at
most countable is the origin. Hence concentration can occur only at
$0$ and at $\infty$.  
 
We mention that action of this type of groups was considered in
\cite{BW} to find nonradial solutions to a Euclidean scalar field
equation. We refer to \cite[\S 1.5]{willemminimax} for a discussion on
the relation between symmetry and compactness in variational problems.

Let us define
\begin{equation}\label{eq:60}
S_{\rm circ}(\l_0)=\inf_{\substack{u\in\Dis\\u\not\equiv 0}}
\frac{{\displaystyle{\int_{\R^N} |\n
      u|^2\,dx-{\l_0}\int_{\R^N}\frac{u^2(x)}{|x|^2}\,dx }}}{{\displaystyle{ 
    \bigg(\int_{\R^N}  |u|^{2^*}dx\bigg)^{2/2^*}}}}.
\end{equation}
The following theorem provides a threshold up to which $J$ satisfies Palais-Smale
condition.

\begin{Theorem}\label{th:ps-circ}
Assume $N\geq 4$ and $\l_0^++\sum_{j=1}^{m}\L_j^+<\frac{(N-2)^2}4$. Let $\{u_n\}_{n\in\N}\subset \Dis$
be a Palais-Smale sequence for
$J$ restricted to $\Dis$, namely 
$$ 
\lim_{n\to\infty}J(u_n)=c<\infty\text{ in }\R\quad\text{and}\quad
\lim_{n\to\infty}J'(u_n)=0\text{ in the dual space }(\Dis)^{\star}. 
$$ 
If 
\begin{align}\label{eq:19-circ}
c<\frac{S_{\rm circ}(\l_0,\L_1,\dots,\L_{m})^{1-\frac{N}{2}}}
{N}
\min\bigg\{S_{\rm circ}(\l_0),
S_{\rm circ}\Big(\l_0+\sum\limits_{\ell=1}^{m}\L_{\ell}\Big)\bigg\}^{\!\!\frac
  N2}\!\!, 
\end{align}
then $\{u_n\}_{n\in\N}$ has a converging subsequence in $\Dis$.
\end{Theorem}
\begin{pf}
Let $\{u_n\}$ be a Palais-Smale sequence for $J$ in $\Dis$, then from
 Sobolev's inequality and (\ref{eq:53}), it is easy to prove that $\{u_n\}$
is a bounded sequence in $\Di$.  Hence, up to
a subsequence, $u_n\rightharpoonup u_0\mbox{ in } \Di$,  
$ u_n\to u_0$  almost everywhere, and, 
from the {\sl Concentration Compactness Principle} by
P. L. Lions \cite{PL1, PL2} 
\begin{align}
\label{eq:21-circ}
&|\nabla u_n|^2\rightharpoonup d\mu\ge |\n
u_0|^2+\mu_0\delta_0+\sum_{j\in {\mathcal
    J} }\mu_{x_j}\delta_{x_j},\\
\label{eq:27-circ}&|u_n|^{2^*}\rightharpoonup d\nu=
|u_0|^{2^*}+\nu_0\delta_0+
+ \sum_{j\in {\mathcal
    J}}\nu_{x_j}\delta_{x_j},\\
\label{eq:33-circ}& \l_0\dfrac{u_n^2}{|x|^2}\rightharpoonup 
d \g_{0}=
\l_{0}\dfrac{u_0^2}{|x|^2}+\g_{0}
\delta_{0}
\end{align}
where  ${\mathcal J}$ is an at most
countable set, $x_j\in\R^N\setminus\{0\}$,
 $\mu_{x_j},\nu_{x_j}\in\R$, $j\in{\mathcal
  J}$, $\mu_0, \nu_0,\g_0\in\R$,  and the above  convergences hold in
the sense of measures. 
We quantify how much the sequence concentrates at infinity  by the  quantities 
$\nu_{\infty}$, $\mu_{\infty}$ defined as in (\ref{eq:59})
and 
$$
\g_{\infty}=\lim_{R\to \infty}\limsup_{n\to
\infty}\int_{|x|>R}\Big(\l_0+\sum_{\ell=1}^m\L_{\ell}\Big)\frac{ u_n^2}{|x|^2}dx.
$$ 
From pointwise convergence of $u_n\in\Dis$ to $u_0$ we deduce that $u_0$ is
invariant by the $\SO(2)\times\SO(N-2)$-action, hence $u_0\in\Dis$.
Moreover, for any $\phi\in C_0(\R^N)$ and for any $\tau\in\SO(2)\times\SO(N-2)$
we have
$$
\int_{\R^N}|u_n|^{2^*}\phi=\int_{\R^N}|u_n|^{2^*}(\phi\circ\tau^{-1})
$$
that, passing to the limit as $n\to\infty$, yields
$$
\sum_{j\in{\mathcal J}}\nu_{x_j}\phi(x_j)=
\sum_{j\in{\mathcal J}}\nu_{x_j}\phi(\tau^{-1}(x_j)).
$$
Arguing as in the proof of Theorem \ref{th:ps}, we deduce that for any 
$j\in{\mathcal J}$, either $\nu_{x_j}=0$ or
 for any $\tau\in \SO(2)\times\SO(N-2)$ there exists
$i\in{\mathcal J}$ such that $\tau(x_j)=x_i$. Namely, if for some $j\in
{\mathcal J}$, $\nu_{x_j}\not=0$, then  ${\mathcal
  O}(x_j)\subset\{x_i:\ i\in {\mathcal J}\}$. When $N\geq 4$, this is not possible
since ${\mathcal J}$ is at most countable whereas  ${\mathcal O}(x_j)$
is more than countable. Therefore $\nu_{x_j}=0$ for all $j\in{\mathcal
  J}$ and we can rewrite (\ref{eq:27-circ}) as 
$$
|u_n|^{2^*}\rightharpoonup d\nu=
|u_0|^{2^*}+\nu_0\delta_0.
$$

\medskip\noindent{\bf Concentration at the origin.}\ We claim that
\begin{equation}\label{eq:20-circ}
\text{either}\quad\nu_0=0\quad\text{or}\quad
\nu_0\geq
\Big(\frac{S_{\rm circ}(\l_0)}{S_{\rm circ}(\l_0,\L_1,\dots,\L_m)}\Big)^{N/2}.
\end{equation}
In order to prove claim (\ref{eq:20-circ}), we consider a smooth 
cut-off function  $\psi_{0}^{\e}\in \Dis$ satisfying $0\leq \psi_{0}^{\e}(x)\leq 1$,
$$
 \psi_{0}^{\e}(x)= 1\quad\mbox{ if }|x|\leq
 \frac{\e}{2},\quad  
\psi_{0}^{\e}(x)= 0\quad\mbox{ if }|x|\geq
\e,\quad\text{and}\quad |\nabla 
\psi_{0}^{\e} |\le \dfrac{4}{\e }. 
$$
From \eqref{eq:60} we obtain
that
\begin{equation*}
\frac{\int_{\R^N} |\n (u_n\psi_{0}^{\e})|^2dx -\l_{0}\int_{\R^N}
{|x|^{-2}}{|\psi_{0}^{\e}|^2u^2_n}\,dx}{\Big(\int_{\R^N}|\psi_{0}^{\e}
u_n|^{2^*}\Big)^{2/2^*}}\geq
S_{\rm circ}(\l_0)
\end{equation*}
hence passing to limit as $n\to\infty$ and $\e\to 0$ we obtain
\begin{equation}\label{eq:34-circ}
\mu_0\geq \g_0+
S_{\rm circ}(\l_0)\big(\nu_0\big)^{2/2^*}.
\end{equation}
On the other hand, testing $J'(u_n)$ with $u_n\psi_{0}^{\e}$ we obtain
\begin{multline*}
\int_{\R^N} |\nabla
  u_n|^2\psi_0^{\e}+\int_{\R^N} u_n \nabla u_n\cdot \nabla 
 \psi_0^{\e} - \l_0\int_{\R^N}\frac{
   u_n^2\psi_0^{\e}}{|x|^2}\\
-\sum_{\ell=1}^m\L_m
\int_{\R^N}\bigg(\media_{S_{r_{\ell}}}\frac{u_n^2(y)\,\psi_0^{\e}(y)}
{|x-y|^2}\,d\sigma(x)\bigg)\,dy
 -S_{\rm
   circ}(\l_0,\L_1,\dots, \L_m)\int_{\R^N}
 \psi_0^{\e}|u_n|^{2^*}=o(1).
\end{multline*}
In view of (\ref{eq:64}), for $\e$ small we have 
\begin{multline*}
\int_{\R^N}\bigg(\media_{S_{r}}\frac{u_n^2(y)\,\psi_0^{\e}(y)}
{|x-y|^2}\,d\sigma(x)\bigg)\,dy=\int_{\R^N}\frac{u_n^2(y)\,\psi_0^{\e}(y)}
{\sqrt{(r^2+|y|^2)^2-4r^2|y'|^2}}\,dx\\\leq \int_{\R^N}\frac{u_n^2(y)\,\psi_0^{\e}(y)}
{|r^2-|y|^2|}\,dx
\leq \frac1{r^2-\e^2}\int_{|y|<\e}u_n^2(y)\,dy\leq{\rm
  const}\,\bigg(\int_{\R^N}|u_n|^{2^*}\bigg)^{2/2^*}
\frac{\e^2}{r^2-\e^2}\ \mathop{\longrightarrow}_{\e\to0}\ 0, 
\end{multline*}
therefore letting $n\to\infty$ and $\e\to 0$  
we infer that
\begin{equation}\label{eq:35-circ}
\mu_{0}-\g_{0}\leq  S_{\rm circ}(\l_0,\L_1,\dots,\L_m)\nu_{0}. 
\end{equation}
From (\ref{eq:34-circ}) and (\ref{eq:35-circ}) we deduce (\ref{eq:20-circ}).

\medskip\noindent{\bf Concentration at infinity.}\ We claim that
\begin{equation}\label{eq:12-circ}
\text{either}\quad\nu_{\infty}=0\quad\text{or}\quad \nu_{\infty}\geq
\bigg(\frac{S_{\rm circ}(\l_0+\sum_{\ell=1}^m
  \L_{\ell})}{S_{\rm circ}(\l_0,\L_1,\L_2,\dots,\L_m)}\bigg)^{{N}/{2}}. 
\end{equation}
Indeed, let $\psi_R$ be a smooth radial cut-off function such that 
$0\le \psi_R(x)\le 1$, $\psi_R(x)=1$ if $|x|>2R$,  $\psi_R(x)=0$ if $|x|<R$,
and $|\nabla \psi_R |\leq 2/R$. 
Using \eqref{eq:56} and taking $\limsup$ as $n\to\infty$ and limit as
$R\to+\infty$, it is easy to show that
\begin{equation}\label{eq:15-circ}
\mu_{\infty}-\g_{\infty}\geq
S_{\rm circ}\Big({\textstyle{\l_0+\sum_{\ell=1}^m\L_{\ell}}}\Big)
\nu_{\infty}^{2/2^*}.
\end{equation} 
On the other hand, testing $J'(u_n)$ with $u_n\psi_R$ we obtain
\begin{align}\label{eq:61}
&\int_{\R^N} |\nabla
  u_n|^2\psi_R+\int_{\R^N} u_n \nabla u_n\cdot \nabla 
 \psi_R - \l_0\int_{\R^N}\frac{
   u_n^2\psi_R}{|x|^2}\\
\nonumber&-\sum_{\ell=1}^m\L_m
\int_{\R^N}\bigg(\media_{S_{r_{\ell}}}\frac{u_n^2(y)\,\psi_R(y)}
{|x-y|^2}\,d\sigma(x)\bigg)\,dy
 -S_{\rm
   circ}(\l_0,\L_1,\dots, \L_m)\int_{\R^N}
 \psi_R|u_n|^{2^*}=o(1).
\end{align}
Let $\bar R>\max\{r_\ell:\ \ell=1,\dots,m\}$. If $R\geq \bar R$, in
view of (\ref{eq:64}) we
have for all $\ell=1.\dots,m$
\begin{gather*}
\bigg|\media_{S_{r_\ell}}\frac{u_n^2(y)\,\psi_R(y)}{|x-y|^2}\,d\sigma(x)-
\frac{u_n^2(y)\,\psi_R(y)}{|y|^2}\bigg|=
\frac{u_n^2(y)\,\psi_R(y)}{|y|^2}
\bigg|\frac{|y|^2-\sqrt{(r^2_{\ell}+|y|^2)^2-4r^2_{\ell}|y'|^2}}{\sqrt{(r^2_{\ell}+|y|^2)^2-4r^2_{\ell}|y'|^2}}\bigg|\\
\leq
\frac{u_n^2(y)\,\psi_R(y)}{|y|^2}\,\frac{r_\ell^2}{|y|^2}\,\frac{6|y|^2+r_{\ell}^2}{|y|^2-r_\ell^2}
\leq \frac{u_n^2(y)\,\psi_R(y)}{|y|^4}\,\frac{7\bar R^4}{\bar R^2-r_{\ell}^2}.
\end{gather*}
Since 
$$
\int_{\R^N}\frac{u_n^2(y)\,\psi_R(y)}{|y|^4}\,dy\leq
\frac1{R^2}\int_{|y|>R}\frac{u_n^2(y)}{|y|^2}\,dy \leq\frac{\rm
  const}{R^2}\ \mathop{\longrightarrow}\limits_{R\to+\infty}\ 0
$$
we deduce that
$$
\int_{\R^N}\media_{S_{r_\ell}}\frac{u_n^2(y)\,\psi_R(y)}{|x-y|^2}\,d\sigma(x)\,dy=
\int_{\R^N}\frac{u_n^2(y)\,\psi_R(y)}{|y|^2}\,dy+o(1)\quad\text{as
}R\to+\infty.
$$
Therefore, letting $n\to\infty$ and $R\to+\infty$ in (\ref{eq:61}), we obtain
\begin{equation}\label{eq:17-circ}
\mu_{\infty}-\g_{\infty}\leq S_{\rm circ}(\l_0,\L_1,\dots,\L_m)\nu_{\infty}. 
\end{equation}
Claim (\ref{eq:12-circ}) follows from (\ref{eq:15-circ}) and (\ref{eq:17-circ}).

\medskip\noindent
As a conclusion we obtain
\begin{align}\label{eq:18-circ} 
c & =  J(u_n)-\frac{1}{2}\langle J'(u_n),u_n\rangle
  +o(1)= \frac{S_{\rm circ}(\l_0,\L_1,\dots,\L_m)}{N}\bigg\{\int_{\R^N}
|u_0|^{2*}dx+\nu_0+\nu_{\infty}\bigg\}.
\end{align}
From (\ref{eq:19-circ}),
(\ref{eq:20-circ}), (\ref{eq:12-circ}),  
and (\ref{eq:18-circ}), we deduce that $\nu_0=0$ and $\nu_{\infty}=0$. Then,
up to a subsequence, $u_n\to u_0$ in $\Dis$.
\end{pf}

\section{Interaction estimates}\label{sec:inter-estim}

\noindent For any $u\in \Di$, let us consider the family of functions
obtained from $u$ by dilation, i.e.
\begin{equation}\label{eq:81}
u_{\mu}(x)=\mu^{-\frac{N-2}2}u(x/\mu), \quad\mu>0.
\end{equation}
The following lemma describes the behavior of
$\int|x+\xi|^{-2}|u^{\lambda}_{\mu}|^2$ as $\mu\to 0$ for any
solution $u^{\lambda}$ of equation (\ref{eq:1}). We mention that
estimates below were obtained in \cite{FT} for radial solutions to
(\ref{eq:1}) (i.e. for functions $w^{(\lambda)}_{\mu}$ in (\ref{eq:ZA})).  
\begin{Lemma}\label{l:interest}
Let $u^{\l}\in\Di$ be a solution to (\ref{eq:1}).
 For any $\xi\in\R^N$ there holds
$$
\int_{\R^N}\frac{|u^{\l}_{\mu}|^2}{|x+\xi|^2}\,dx=
\begin{cases}
{ {\frac{\mu^2}{|\xi|^2}\int_{\R^N}|u^{\l}|^2\,dx
+o\left(\mu^2\right)}}& \text{if }\l<\frac{N(N-4)}4
,\\[10pt]
\kappa_{\infty}(u^{\lambda})^2\, { {\frac{\mu^2|\ln\mu|}{|\xi|^2}}}
 +o(\mu^2|\ln\mu|)
& \text{if }\l=\frac{N(N-4)}4,\\[10pt]
\kappa_{\infty}(u^{\lambda})^2\,\beta_{\l,N}\,
\frac{\mu^{\sqrt{(N-2)^2-4\lambda}}}
{|\xi|^{\sqrt{(N-2)^2-4\lambda}}}+o\big(\mu^{\sqrt{(N-2)^2-4\lambda}}\big)
 &
\text{if }\l>\frac{N(N-4)}4,
\end{cases}
$$
as $\mu\to 0$,
where $$
\beta_{\l,N}=\int_{\R^N}\frac{dx}
{|x|^2|x-e_1|^{N-2+\sqrt{(N-2)^2-4\l}}}, \quad e_1=(1,0,\dots,0)\in\R^N
.$$
\end{Lemma}

\begin{pf}
We have that
\begin{align}\label{eq:7}
\int_{\R^N}\frac{|u_{\mu}^{\l}|^2}{|x+\xi|^2}\,dx
&=\mu^2\int_{\R^N}\frac{|u^{\l}|^2}{|\mu x+\xi|^2}\,dx
 =\mu^2\int_{|x|<\frac{|\xi|}{2\mu}}\frac{|u^{\l}|^2}{|\mu x+\xi|^2}\,dx
+\mu^2\int_{|x|>\frac{|\xi|}{2\mu}}\frac{|u^{\l}|^2}{|\mu x+\xi|^2}\,dx.
\end{align}
For $\l<\frac{N(N-4)}4$, from (\ref{eq:2}) we have that
 $ u^{\l}\in L^2(\R^N)$. 
From (\ref{eq:2}) and \cite[Proof of Lemma 3.4]{FT} we deduce
\begin{align*}
\bigg|\int_{|x|<\frac{|\xi|}{2\mu}}&|u^{\l}(x)|^2\bigg[\frac{1}{|\mu
  x+\xi|^2}-\frac{1}{|\xi|^2}\bigg]\,dx\bigg|
\\
&\leq
\kappa(u^{\lambda})^2\int_{|x|<\frac{|\xi|}{2\mu}}|w^{\l}_1(x)|^2
\bigg|\frac{1}{|\mu  x+\xi|^2}-\frac{1}{|\xi|^2}\bigg|\,dx=o(1)
\end{align*}
as $\mu\to 0$  and hence, since $u^{\l}\in L^2(\R^N)$,
\begin{align}\label{eq:23}
\int_{|x|<\frac{|\xi|}{2\mu}}
|u^{\l}(x)|^2\frac{dx}{|\mu
x+\xi|^2}=
\frac{1}{|\xi|^2}\int_{|x|<\frac{|\xi|}{2\mu}}|u^{\l}(x)|^2
\,dx+o(1)=\frac{1}{|\xi|^2}\int_{\R^N}|u^{\l}(x)|^2
\,dx
+o(1).\quad
\end{align}
On the other hand, from \cite[Proof of Lemma 3.4]{FT}  we have 
\begin{align}\label{eq:24}
\mu^2\int_{|x|>\frac{|\xi|}{2\mu}}\frac{|u^{\l}|^2}{|\mu x+\xi|^2}\,dx
&\leq \kappa(u^{\lambda})^2\mu^{2-N}\int_{|x-\xi|\geq\frac{|\xi|}2}
\bigg|w^{\lambda}_1\Big(\frac{x-\xi}\mu\Big)\bigg|^2\,\frac{dx}{|x|^2}\\
&\nonumber
=O\big(\mu^{\nu_{\lambda}(N-2)}\big)=o(\mu^2).
\end{align}
From (\ref{eq:7}), (\ref{eq:23}), and (\ref{eq:24}) we deduce that
\[
\int_{\R^N}\frac{|u^{\l}|^2}{|x+\xi|^2}\,dx
=\frac{\mu^2}{|\xi|^2}\int_{\R^N}|u^{\l}(x)|^2
\,dx
 +o(\mu^2).
\]
For $\l=\frac{N(N-4)}4$, from (\ref{eq:2}) and  \cite[Proof of Lemma
3.4]{FT} 
we deduce that
\begin{align}\label{eq:26}
\int_{|x|<\frac{|\xi|}{2\mu}}
|u^{\l}(x)|^2\frac{dx}{|\mu
x+\xi|^2}=
\frac{1}{|\xi|^2}\int_{|x|<\frac{|\xi|}{2\mu}}|u^{\l}(x)|^2
\,dx+O(1).
\end{align}
 On the other hand, from (\ref{eq:4}) we obtain 
\begin{align}\label{eq:29}
\int_{|x|<\frac{|\xi|}{2\mu}}
|u^{\l}(x)|^2
\,dx&=\int_{1<|x|<\frac{|\xi|}{2\mu}}
|u^{\l}(x)|^2
\,dx+\int_{|x|<1}
|u^{\l}(x)|^2
\,dx\\
&\nonumber
=\kappa_{\infty}(u^{\lambda})^2\int_{1<|x|<\frac{|\xi|}{2\mu}}|x|^{-N}\,dx
+O\bigg( \int_{1<|x|<\frac{|\xi|}{2\mu}}|x|^{-N-2\a}\,dx\bigg)+O(1)
\\
&\nonumber
=\kappa_{\infty}(u^{\lambda})^2|\ln\mu|+O(1).
\end{align}
Arguing as above (see (\ref{eq:24})), we obtain
\begin{align}\label{eq:31}
\mu^2\int_{|x|>\frac{|\xi|}{2\mu}}\frac{|u^{\l}|^2}{|\mu x+\xi|^2}\,dx
=O(\mu^2).
\end{align}
Gathering  (\ref{eq:7}), (\ref{eq:26}), (\ref{eq:29}) and
(\ref{eq:31}) we deduce that 
\[
\int_{\R^N}\frac{|u_{\mu}^{\l}|^2}{|x+\xi|^2}\,dx
=\kappa_{\infty}(u^{\lambda})^2\, \frac{\mu^2|\ln\mu|}{|\xi|^2}
 +o(\mu^2|\ln\mu|).
\]
For $\l>\frac{N(N-4)}4$, in view of (\ref{eq:4}) we have that
\begin{align}\label{eq:32}
\int_{\R^N}&\frac{|u_{\mu}^{\l}|^2}{|x+\xi|^2}\,dx
=\mu^{\nu_{\l}(N-2)}\bigg[
\kappa_{\infty}^2(u^{\lambda})\int_{\R^N}
\frac{1}
{|x|^2|x-\xi|^{(N-2)(1+\nu_{\l})}}+o(1)
\bigg]\,dx
\end{align}
As observed in  \cite[Proof of Lemma
3.4]{FT}, the function 
$$
\varphi(\xi):=\int_{\R^N}\frac{dx}
{|x|^2|x-\xi|^{(N-2)(1+\nu_{\l})}}
$$
can be written as 
\begin{equation}\label{eq:40}
\varphi(\xi)=|\xi|^{-\sqrt{(N-2)^2-4\lambda}}\varphi(\xi/|\xi|)=
|\xi|^{-\sqrt{(N-2)^2-4\lambda}}
\varphi(e_1).
\end{equation}
(\ref{eq:32}) and (\ref{eq:40}) yield the required estimate for
$\l>\frac{N(N-4)}4$. 
\end{pf}

\noindent Let us now study the interaction between two minimizers of
(\ref{eq:minimiz}), 
i.e. functions $z_{\mu}^{\l}$ in (\ref{eq:44}), centered at different
points as $\mu\to 0$. To this aim we note that a direct calculation
yields
\begin{align}\label{eq:8}
&z_1^{\l}(x)=|x|^{-\frac{N-2}2(1+\nu_{\l})}\big[\a_{\l,N}
+O(|x|^{-\a})\big],\\
\label{eq:14}
&\n z_1^{\l}(x)=|x|^{-\frac{N+2}2-\nu_{\l}\frac{N-2}2}x\big[-\a_{\l,N}
{\textstyle{\frac{N-2}2}}(1+\nu_{\l})
+O(|x|^{-\a})\big],
\end{align}
for all $0<\a\leq 2\nu_{\l}$. From (\ref{eq:44}), (\ref{eq:8}), and
(\ref{eq:14}), it is easy to deduce the following result.
\begin{Lemma}\label{l:gradgrad}
 For any $\l\in\big(0,(N-2)^2/4\big)$ and $\xi,\zeta\in\R^N$,
 $\xi\not=0$,  there holds
$$
\int_{\R^N} \frac{z_{\mu}^{\l}(x)z_{\mu}^{\l}(x+\xi)}{|x+\zeta|^2}\,dx
=\mu^{\sqrt{(N-2)^2-4\lambda}}\bigg[\a_{\l,N}^2
\int_{\R^N}\frac{dx}{|x|^{(1+\nu_{\l})\frac{N-2}2}
|x+\xi|^{(1+\nu_{\l})\frac{N-2}2}|x+\zeta|^2}+o(1)\bigg]
$$
and
\begin{align*}
&\int_{\R^N} \n z_{\mu}^{\l}(x)\cdot  \n z_{\mu}^{\l}(x+\xi)\,dx
= 
\mu^{\sqrt{(N-2)^2-4\lambda}}\Big[
\a_{\l,N}^2{\textstyle{\frac{(N-2)^2}4}}(1+\nu_{\l})^2\gamma_{\l,N}
|\xi|^{-\sqrt{(N-2)^2-4\lambda}}+o(1)\Big]
\end{align*}
as $\mu\to 0$,
where $$
\gamma_{\l,N}=\int_{\R^N}\frac{x\cdot (x+e_1)\,dx}
{|x|^{\frac{N+2}2+\nu_{\l}\frac{N-2}2}|x-e_1|^{\frac{N+2}2+\nu_{\l}\frac{N-2}2}}, \quad e_1=(1,0,\dots,0)\in\R^N
.$$
\end{Lemma}

\section{Comparison between concentration levels for $J$ and proof of
  Theorem \ref{t:achcirc}}\label{sec:comp-betw-conc}  

In order to compare the level $S_{\rm circ}(\l_0,\L_1,\dots,\L_m)$
with the level $S_{\rm circ}(\l_0)$ of possible concentration at $0$,
we need the following lemma, which states that the infimum in
\eqref{eq:60} is achieved if $N\geq4$. Such a result does not come
unexpected, since it can be seen as the analogue of Lemma~\ref{l:ska}
when $k=\infty$; indeed when $k$ becomes larger and larger, assumption
  $S_k(\l)<k^{2/N}S$ of Lemma \ref{l:ska} is weakened till it is  no
  more needed in the limiting problem corresponding to singularities
  spread over circles. 

\begin{Lemma}\label{l:scircach}
For any $\l_0\in(-\infty,(N-2)^2/4)$ and $N\geq 4$, the infimum in \eqref{eq:60} is
achieved.
\end{Lemma}
\begin{pf}
Hardy's and Sobolev's inequalities imply that
$S_{\rm circ}(\l_0)\geq\big (1-\frac{4\l_0}{(N-2)^2}\big)S>0$. 
Let $\{u_n\}_n\subset\Dis$ be a minimizing sequence such that
$\int_{\R^N}|u_n|^{2^*}=1$. 
By virtue of the Ekeland's variational
principle we can assume that $\{u_n\}_n$ is a Palais-Smale
sequence for the functional 
$$
F(u)=\frac12\int_{\R^N}|\n
u|^2\,dx-\frac{\l_0}2\int_{\R^N}\frac{|u(x)|^2}{|x|^2}\,dx
-\frac{S_{\rm circ}(\l_0)}{2^*}\int_{\R^N}|u(x)|^{2^*}\,dx,\quad u\in\Dis,
$$
i.e. 
$$ 
\lim_{n\to\infty}F(u_n)=\frac{S_{\rm circ}(\l_0)}{N}
\text{ in }\R\quad\text{and}\quad
\lim_{n\to\infty}F'(u_n)=0\text{ in the dual space }(\Dis)^{\star}. 
$$ 
Let 
$$
w_n(x)=\sigma_n^{-\frac{N-2}2}u_n(\sigma_n^{-1}x)
$$
where $\sigma_n$ is chosen in such a way that
\begin{equation}\label{eq:68}
\int_{B(0,1)}|w_n(x)|^{2^*}\,dx=\int_{B(0,\sigma_n^{-1})}
|u_n(x)|^{2^*}\,dx=\frac12.
\end{equation}
Scaling invariance ensures that $\{w_n\}_n\subset\Dis$ is also a
minimizing sequence for \eqref{eq:60},
\begin{equation}\label{eq:72}
\int_{\R^N}|w_n(x)|^2\,dx=1,\quad
\lim_{n\to+\infty}\int_{\R^N}\bigg(|\nabla
w_n(x)|^2-\l_0\frac{|w_n(x)|^2}{|x|^2}\bigg)\,dx =S_{\rm circ}(\l_0),
\end{equation}
and
\begin{equation}\label{eq:65}
F'(w_n)\to 0 \text{ in the dual space }(\Dis)^{\star}.
\end{equation}
Since $\{w_n\}_n$
is bounded in $\Di$, up to
a subsequence $w_n$ converges to $w$ weakly in $\Di$ and
almost everywhere. Pointwise convergence implies that
$w\in\Dis$. From the {\sl Concentration Compactness Principle} by
P. L. Lions \cite{PL1, PL2} and taking into account that, when $N\geq 4$,   
 $\SO(2)\times\SO(N-2)$-invariant functions
 can concentrate only at $0$ and at $\infty$, as already
pointed out in the proof of Theorem \ref{th:ps-circ}, we have
\begin{align}\label{eq:63}
|\nabla w_n|^2\rightharpoonup d\mu\ge |\n
w|^2+\mu_0\delta_0,\ \
|w_n|^{2^*}\rightharpoonup d\nu=
|w|^{2^*}+\nu_0\delta_0,\ \ \text{and}\  
 \dfrac{w_n^2}{|x|^2}\rightharpoonup 
d \g=
\dfrac{w^2}{|x|^2}+\g_{0}
\delta_{0}.
\end{align}
The amount of concentration at infinity is  quantified
by the following numbers 
\begin{equation*}
\nu_{\infty}=\lim_{R\to \infty}\limsup_{n\to
\infty}\int_{|x|>R}|w_n|^{2^*}dx,\quad \mu_{\infty}=\lim_{R\to
\infty}\limsup_{n\to \infty}\int_{|x|>R}|\n w_n|^{2}dx,
\end{equation*}
and  $$
\g_{\infty}=\lim_{R\to \infty}\limsup_{n\to
\infty}\int_{|x|>R}\frac{ w_n^2}{|x|^2}dx.
$$
Arguing as we did to prove (\ref{eq:34-circ}) and (\ref{eq:15-circ}),
we can easily obtain 
\begin{equation}\label{eq:76}
\nu_0^{\frac2N}(\mu_0-\l_0\gamma_0)\geq S_{\rm circ}(\l_0)\nu_0
\end{equation}
and 
\begin{equation}\label{eq:80}
(\mu_{\infty}-\l_0\gamma_{\infty})\geq S_{\rm
  circ}(\l_0)\nu_{\infty}^{\frac{2}{2^*}}. 
\end{equation}

\medskip\noindent
{\bf Step 1: \rm we prove that $w\not\equiv 0$.} 
  By contradiction, assume that $w\equiv 0$. Then from
(\ref{eq:72}--\ref{eq:63})  we deduce 
\begin{equation}\label{eq:73}
1=\nu_0+\nu_{\infty}
\end{equation}
and
\begin{equation}\label{eq:74}
S_{\rm circ}(\l_0)=\int_{\R^N} d\mu-\l_0\int_{\R^N} d\gamma
+\mu_{\infty}-\l_0\gamma_{\infty}.
\end{equation}
Let $\psi_R$ be a smooth radial cut-off function such that 
$0\le \psi_R(x)\le 1$, $\psi_R(x)=1$ if $|x|>2R$,  $\psi_R(x)=0$ if $|x|<R$,
and $|\nabla \psi_R |\leq 2/R$. From (\ref{eq:65}), we have
\begin{align*}
o(1)&=\big\langle F'(w_n),w_n \psi_R\big\rangle\\
&=\int_{\R^N}|\nabla
w_n|^2\psi_R+\int_{\R^N}w_n\n w_n\cdot\n
\psi_R-\l_0\int_{\R^N}\frac{w_n \psi_R}{|x|^2}
-S_{\rm circ}(\l_0)\int_{\R^N}|w_n|^{2^*}\psi_R.
\end{align*}
Taking $\limsup$ as $n\to\infty$ and limit as
$R\to+\infty$, we find
\begin{equation}\label{eq:75}
\mu_{\infty}-\l_0\gamma_{\infty}=S_{\rm circ}(\l_0)\nu_{\infty}.
\end{equation}
From (\ref{eq:74}--\ref{eq:75}) it follows that 
\begin{equation}\label{eq:77}
\int_{\R^N} d\mu-\l_0\int_{\R^N} d\gamma
=S_{\rm circ}(\l_0)(1-\nu_{\infty}).
\end{equation}
From (\ref{eq:73}), (\ref{eq:76}), and (\ref{eq:77}), it follows that
\begin{align}\label{eq:78}
1-\nu_{\infty}&=\nu_0\leq S_{\rm circ}(\l_0)^{-1}
\nu_0^{\frac2N}(\mu_0-\l_0\gamma_0)\leq  S_{\rm circ}(\l_0)^{-1}
\nu_0^{\frac2N}\bigg(\int_{\R^N} d\mu-\l_0\int_{\R^N}
d\gamma\bigg)\\
\nonumber&=\nu_0^{\frac2N}(1-\nu_{\infty})=(1-\nu_{\infty})^{1+\frac2N}.
\end{align}
On the other
hand, from \eqref{eq:68} we have
$$
\int_{\R^N\setminus\{B(0,R)\}}|w_n|^{2^*}\,dx=1-
\int_{B(0,R)}|w_n|^{2^*}\,dx\leq 1-
\int_{B(0,1)}|w_n|^{2^*}\,dx=\frac12\quad\text{for all }R>1,
$$
hence 
\begin{equation}\label{eq:79}
\nu_{\infty}\leq \frac12.
\end{equation}
From (\ref{eq:78}--\ref{eq:79}) we deduce that
$\nu_{\infty}=0$. From (\ref{eq:73}) it follows that $\nu_0=1$. Therefore
\eqref{eq:68} implies
$$
\frac12=\int_{B(0,1)}|w_n|^{2^*}\
\mathop{\longrightarrow}\limits_{n\to\infty}\ d\nu(B(0,1))=1,
$$
thus giving rise to a contradiction.
 
\medskip\noindent
{\bf Step 2:}  we prove that 
$$
\int_{\R^N}|w|^{2^*}=1\quad\text{and}\quad
\int_{\R^N}|\n
w|^{2}\,dx-\l_0\int_{\R^N}\frac{|w|^2}{|x|^2}\,dx=
S_{\rm circ}(\l_0).
$$ 
Let $\rho=\int_{\R^N}|w|^{2^*}$. From (\ref{eq:63}), we have that
$1=\rho+\nu_0+\nu_{\infty}$ and, in view of step 1, $\rho\in(0,1]$,
i.e. $\nu_0+\nu_{\infty}\in[0,1)$.  Since
\begin{align*}
S_{\rm circ}(\l_0)&=\int_{\R^N} d\mu-\l_0\int_{\R^N} d\gamma
+\mu_{\infty}-\l_0\gamma_{\infty}\\
&\geq \int_{\R^N}|\n
w|^{2}\,dx-\l_0\int_{\R^N}\frac{|w|^2}{|x|^2}\,dx+\mu_0-\l_0\gamma_0+\mu_{\infty}-\l_0\gamma_{\infty},
\end{align*}
from (\ref{eq:76}), (\ref{eq:80}), and concavity of the function
$t\mapsto t^{2/2^*}$ we deduce
\begin{align*}
\int_{\R^N}|\n
w|^{2}\,dx-\l_0\int_{\R^N}\frac{|w|^2}{|x|^2}\,dx&\leq S_{\rm
  circ}(\l_0)\big(
1-\nu_0^{\frac2{2^*}}-\nu_{\infty}^{\frac2{2^*}}\big)\leq 
S_{\rm circ}(\l_0)\big(
1-\nu_0-\nu_{\infty}\big)^{\frac2{2^*}}\\
&=S_{\rm circ}(\l_0)\rho^{2/2^*}\leq \int_{\R^N}|\n
w|^{2}\,dx-\l_0\int_{\R^N}\frac{|w|^2}{|x|^2}\,dx.
\end{align*}
Hence all the above inequalities are indeed equalities; in particular 
$\big(
1-\nu_0^{\frac2{2^*}}-\nu_{\infty}^{\frac2{2^*}}\big)=\big(
1-\nu_0-\nu_{\infty}\big)^{\frac2{2^*}}$ which is possible only when
$\nu_0=\nu_{\infty}=0$. 
Therefore $\rho=1$ and 
$$
\int_{\R^N}|\n
w|^{2}\,dx-\l_0\int_{\R^N}\frac{|w|^2}{|x|^2}\,dx=
S_{\rm circ}(\l_0),
$$
i.e. $w$ attains the infimum.
\end{pf}

\noindent We now provide a sufficient condition for the
infimum 
in (\ref{eq:56}) to stay below the level $S_{\rm circ}(\l_0)$, at
which possible concentration at $0$ can occur.

\begin{Lemma}\label{l:comp-betw-conc}
Let $\l_0,\Lambda_1,\dots,\Lambda_m\in\R$, $r_1,r_2,\dots,r_m\in\R^+$
satisfy \eqref{eq:58}. Then  
$$
S_{\rm circ}(\l_0,\L_1,\dots,\L_m)<S_{\rm circ}(\l_0).
$$
\end{Lemma}
\begin{pf}
From Lemma \ref{l:scircach}, we have that $S_{\rm circ}(\l_0)$ is attained by some
$u^{\l_0}\in\Dis$. By homogeneity of the Rayleigh quotient, we can
assume $\int |u^{\l_0}|^{2^*}=1$. Moreover, the function
  $v^{\l_0}=S_{\rm circ}(\l_0)^{1/(2^*-2)}|u^{\l_0}|$ is a
  nonnegative solution to (\ref{eq:1}), hence we can apply
  Lemma~\ref{l:interest} to study the behavior of 
  $\int_{\R^N}\frac{|u^{\l_0}_{\mu}|^2}{|x+\xi|^2}\,dx$ as $\mu\to
  0$, where $u^{\l_0}_{\mu}$ are defined in (\ref{eq:81}). Hence for
  some positive constant $\tilde\kappa$
\begin{multline*}
\int_{\R^N}\bigg(\media_{S_{r_\ell}}\frac{|u_{\mu}^{\l_0}(y)|^2}
{|x-y|^2}\,d\sigma(x)\bigg)\,dy\\[10pt]
=
\begin{cases}
{ {\frac{\mu^2}{r_{\ell}^2}\int_{\R^N}|u_{1}^{\l_0}|^2\,dx
+o\left(\mu^2\right)}}& \text{if }\l_0<\frac{N(N-4)}4
,\\[10pt]
\tilde\kappa^2\, { {\frac{\mu^2|\ln\mu|}{r_{\ell}^2}}}
 +o(\mu^2|\ln\mu|)
& \text{if }\l_0=\frac{N(N-4)}4,\\[10pt]
\tilde\kappa^2\,\beta_{\l_0,N}\,
\mu^{\sqrt{(N-2)^2-4\lambda_0}}
|r_{\ell}|^{-\sqrt{(N-2)^2-4\lambda_0}}+o\big(\mu^{\sqrt{(N-2)^2-4\lambda_0}}\big)
 &
\text{if }\l_0>\frac{N(N-4)}4,
\end{cases}
\end{multline*}
as $\mu\to 0$,
where $\beta_{\l_0,N}$ is defined in Lemma \ref{l:interest}.
Therefore
\begin{multline}\label{eq:70}
S_{\rm circ}(\l_0,\L_1,\dots,\L_m)\\
\leq 
\int_{\R^N} |\n
      u_{\mu}^{\l_0}|^2\,dy-{\l_0}\int_{\R^N}\frac{|
        u_{\mu}^{\l_0}(y)|^2}{|y|^2}\,dy 
 -\sum_{\ell=1}^m
 {\L_{\ell}}\int_{\R^N}\bigg(\media_{S_{r_{\ell}}}\frac{| u_{\mu}^{\l_0} (y)|^2}
{|x-y|^2}\,d\sigma(x)\bigg)\,dy\\
= S_{\rm circ}(\l_0)-
\begin{cases}
\mu^2\big(\int_{\R^N}|u_{1}^{\l_0}|^2\big)\bigg({\displaystyle{\sum_{\ell=1}^m
 \frac{\L_{\ell}}{r_{\ell}^2}+o(1)}}\bigg)& \text{if }\l_0<\frac{N(N-4)}4
,\\[10pt]
\mu^2|\ln\mu|\tilde\kappa^2\, \bigg({\displaystyle{\sum_{\ell=1}^m
 \frac{\L_{\ell}}{r_{\ell}^2}+o(1)}}\bigg)
& \text{if }\l_0=\frac{N(N-4)}4,\\[10pt]
\mu^{\sqrt{(N-2)^2-4\lambda_0}}\tilde\kappa^2\,\beta_{\l_0,N}\,
\bigg({\displaystyle{\sum_{\ell=1}^m
 \frac{\L_{\ell}}{|r_{\ell}|^{\sqrt{(N-2)^2-4\lambda_0}}}+o(1)}}\bigg)
 &
\text{if }\l_0>\frac{N(N-4)}4,
\end{cases}
\end{multline}
as $\mu\to 0$. Taking $\mu$ sufficiently small, assumption
\eqref{eq:58} yields $S_{\rm circ}(\l_0,\L_1,\dots,\L_m)< S_{\rm circ}(\l_0)$.
~\end{pf}

\begin{pfn}{Theorem \ref{t:achcirc}}
 Let $\{u_n\}_n\subset\Dis$ be a minimizing sequence for
\eqref{eq:56}. From the homogeneity of the quotient there is no
restriction requiring
$\|u_n\|_{L^{2^*}(\R^N)}=1$. Moreover from Ekeland's variational
principle we can assume that $\{u_n\}_n\subset\Dis$ is a Palais-Smale
sequence, more precisely  $J'(u_n)\to 0$ in $(\Dis)^{\star}$ and 
$J(u_n)\to \frac1N S_{\rm circ}(\l_0\L_1,\dots, \L_m)$. 
From assumption \eqref{eq:57} we deduce that
\begin{equation}\label{eq:67}
S_{\rm circ}(\l_0)\leq S_{\rm
  circ}\Big(\l_0+\sum\limits_{\ell=1}^{m}\L_{\ell}\Big).
\end{equation}
From  Lemma \ref{l:comp-betw-conc} and  (\ref{eq:67}), it follows that
the level of the minimizing Palais-Smale sequence satisfies assumption
(\ref{eq:19-circ}). Hence from Theorem \ref{th:ps-circ}, 
$\{u_n\}_{n\in\N}$ has a  subsequence strongly
converging to some $u_0\in\Dis$ 
such that 
$J(u_0)=\frac1N S_{\rm circ}(\l_0,\L_1,\dots,\L_m)$.
Hence $u_0$ achieves the infimum in \eqref{eq:56}. Since
$J$ is even, also $|u_0|$ is a minimizer in \eqref{eq:56}
and then $v_0=S_{\rm circ}(\l_0,\L_1,\dots,\L_m)^{1/(2^*-2)}|u_0|$ is a
nonnegative solution to equation \eqref{eq:55}. The maximum principle
implies the positivity  outside singular circles of such a solution.
~\end{pfn}

\section{Limit of $S_k(\l)$ as $k\to\infty$.}\label{sec:limit-s_kl-as}

\noindent Since $\Dis\subset \Dik$, for any $\l\in(-\infty,(n-2)^2/4)$,  there holds
\begin{equation}\label{eq:ven1}
0<S(\l)\leq S_k(\l)\leq S_{\rm circ}(\l).
\end{equation}
\noindent From \eqref{eq:ven1} and Lemma \ref{l:ska}, it follows easily the following result.

\begin{Lemma}\label{l:ven2}
Let $\l\in(-\infty,(n-2)^2/4)$ and $N\geq 4$. Then there exists $\bar k=\bar k(\l,N)$ such that $S_k(\l)$ is achieved for all $k\geq\bar k$.
\end{Lemma}

Let us now study the limit of $S_k(\l)$ as $k\to\infty$. Theorem \ref{t:ven} provides convergence of $S_k(\l)$ to $S_{\rm circ}(\l)$. To prove it we will need the following proposition.

\begin{Proposition}\label{p:sab}
Let $\l\in(-\infty,(n-2)^2/4)$ and let $\{w_k\}_k$ be a sequence in $\Di$  such that $w_k\in\Dik$, 
\begin{equation}\label{eq:pom1}
\int_{\R^N}|w_k|^{2^*}=1,\quad Q_{\l}(w_k)=S_k(\l),
\end{equation}
and $w_k$ converges weakly to $0$ in $\Di$ (at least along a subsequence). Then, for any $r>0$ and $\e\in(-r,r)$, there exists $\rho$ such that $0<|\rho|<|\e|$ and , for a subsequence,
\begin{align*}
&\text{either}\quad \int_{B(0,r+\rho)}|\n w_k|^2\to 0,\quad
\int_{B(0,r+\rho)}|w_k|^{2^*}\to 0,\quad\text{and }\int_{B(0,r+\rho)}\frac{|w_k|^2}{|x|^2}\to 0,\\
&
\text{or}\quad 
\int_{\R^N\setminus B(0,r+\rho)}|\n w_k|^2\to 0,\quad
\int_{\R^N\setminus B(0,r+\rho)}|w_k|^{2^*}\to 0,\quad\text{and }\int_{\R^N\setminus B(0,r+\rho)}\frac{|w_k|^2}{|x|^2}\to 0.
\end{align*}
\end{Proposition}
\begin{pf}
An analogous result is proved in \cite{terracini} for minimizing sequences of quotient \eqref{eq:minimiz}. Since the proof of Proposition \ref{p:sab} is similar, we will be sketchy.

Let $\e\in(0,r)$ (the proof for $\e$ negative is similar). Since
$$
\int_r^{r+\e}d\rho\int_{\partial B(0,\rho)}|\n w_k|^2 =\int_{B(0,r+\e)\setminus B(0,r)}
|\n w_k|^2,
$$
we can choose $\rho\in(0,\e)$ such that, for infinitely many $k$'s (i.e. along a subsequence still denoted as $\{w_k\}_k$) 
\begin{equation}\label{eq:po1}
\int_{\partial B(0,r+\rho)}|\n w_k|^2\leq \frac2{\e} 
\int_{B(0,r+\e)\setminus B(0,r)}
|\n w_k|^2.
\end{equation}
From the uniform bound of $S_k(\l)$ (see \eqref{eq:ven1}) and equivalence of $Q_{\l}$ to $\Di$-norm, it follows that 
\begin{equation}\label{eq:po2}
\int_{B(0,r+\e)\setminus B(0,r)}
|\n w_k|^2\leq {\rm const}\,\, Q_{\l}(w_k)=
{\rm const}\,\, S_k(\l)\leq {\rm const}.
\end{equation}
From (\ref{eq:po1}--\ref{eq:po2}), it follows that $\big\{w_k\big|_{\partial B(0,r+\rho)}\big\}_k$ is bounded in $H^1\big(\partial B(0,r+\rho)\big)$. By compactness of the embedding $H^1\big(\partial B(0,r+\rho)\big)\hookrightarrow H^{1/2}\big(\partial B(0,r+\rho)\big)$
 and weak convergence to $0$, we conclude that, up to subsequence,  
$\big\{w_k\big|_{\partial B(0,r+\rho)}\big\}_k$  converges strongly to $0$ in $H^{1/2}\big(\partial B(0,r+\rho)\big)$. 

Let $w^1_k$ (respectively $w^2_k$) be the harmonic functions in $B(0,r+\e)\setminus B(0,r+\rho)$ (respectively $B(0,r+\rho)\setminus B(0,r-\e)$) such that
$$
\begin{cases}
w^1_k=w_k&\text{on }\partial B(0,r+\rho),\\
w^1_k=0&\text{on }\partial B(0,r+\e),
\end{cases}
\qquad
\begin{cases}
w^2_k=w_k&\text{on }\partial B(0,r+\rho),\\
w^2_k=0&\text{on }\partial B(0,r-\e).
\end{cases}
$$
By continuity of $\D^{-1}$ we have that $w^1_k\to 0$ strongly in  $H^1(B(0,r+\e)\setminus B(0,r+\rho))$ and $w^2_k\to 0$ strongly in  $H^1(B(0,r+\rho)\setminus B(0,r-\e))$. Moreover symmetry properties of boundary data ensure that $w^1_k,w^2_k \in\Dik$. Let us now set
$$
\begin{cases}
u^1_k=w_k&\text{in }B(0,r+\rho),\\
u^1_k=w^1_k&\text{in }B(0,r+\e)\setminus B(0,r+\rho),\\
u^1_k=0&\text{in }\R^N\setminus  B(0,r+\e),
\end{cases}
\qquad
\begin{cases}
u^2_k=w_k&\text{in }\R^N\setminus B(0,r+\rho),\\
u^2_k=w^2_k&\text{in }B(0,r+\rho)\setminus B(0,r-\e),\\
u^2_k=0&\text{in } B(0,r-\e).
\end{cases}
$$
Direct computations yield
\begin{align*}
&Q_{\l}(u^1_k)=S_k(\l)\int_{B(0,r+\rho)}|w_k|^{2^*}\,dx+o(1),\\
&Q_{\l}(u^2_k)=S_k(\l)\int_{\R^N\setminus B(0,r+\rho)}|w_k|^{2^*}\,dx+o(1),\\
&Q_{\l}(u^1_k)+Q_{\l}(u^2_k)=S_k(\l)\int_{\R^N}|w_k|^{2^*}\,dx+o(1)=
S_k(\l)+o(1)=Q_{\l}(w_k)+o(1),\\
&
\int_{\R^N}|u^1_k|^{2^*}\,dx+\int_{\R^N}|u^2_k|^{2^*}\,dx=\int_{\R^N}|w_k|^{2^*}\,dx+o(1).
\end{align*}
We claim that either $Q_{\l}(u^1_k)\to 0$ or $Q_{\l}(u^2_k)\to 0$
along some subsequence. Indeed, assume that $Q_{\l}(u^1_k)\not\to 0$
along any subsequence, i.e. $Q_{\l}(u^1_k)$ stays bounded below away from $0$.  
From above, \eqref{eq:pom1} and \eqref{eq:minimizk}, it follows
\begin{align*}
\frac{Q_{\l}(u^1_k)}{\|u^1_k\|_{2^*}^2}&=\frac{Q_{\l}(w_k)-Q_{\l}(u^2_k)+o(1) }{\big(\|w_k\|_{2^*}^{2^*}-\|u^2_k\|_{2^*}^{2^*}+o(1)\big)^{2/2^*}  } \leq S_k(\l)
\frac{Q_{\l}(w_k)-Q_{\l}(u^2_k)+o(1) }{\big(Q_{\l}(w_k)^{2^*/2}-Q_{\l}(u^2_k)^{2^*/2}+o(1)\big)^{2/2^*}}<S_k(\l)
\end{align*}
in contradiction with \eqref{eq:minimizk}, unless $Q_{\l}(u^2_k)\to 0$ along some subsequence. 
The claim is thereby proved. The statement of the proposition follows
from equivalence to norm of $Q_{\l}$, Hardy's and Sobolev's inequalities.
\end{pf}

\begin{Theorem}\label{t:ven}
Let $\l\in(-\infty,(n-2)^2/4)$ and $N\geq 4$. Then $\lim_{k\to+\infty}S_k(\l)=S_{\rm circ}(\l)$.
\end{Theorem}

\begin{pf}
From Lemma \ref{l:ven2} we know that, for $k$ sufficiently large, $S_k(\l)$ is achieved, hence there exists some $u_k\in\Dik$ such that
\begin{equation*}
\int_{\R^N}|u_k|^{2^*}=1\quad\text{and}\quad Q_{\l}(u_k)=S_k(\l).
\end{equation*}
From the uniform bound of $S_k(\l)$ (see \eqref{eq:ven1}) and equivalence of $Q_{\l}$ to $\Di$-norm, it follows that $\{u_k\}_k$ is bounded in $\Di$. Let us set 
\begin{equation}\label{eq:ser1}
\tilde u_k(x)=R_k^{-\frac{N-2}2}u_k\Big(\frac x{R_k}\Big)\quad\text{and}\quad v_k(x)=(S_k(\l))^{\frac1{2^*-2}}\tilde u_k(x)
\end{equation}
where $R_k$ is chosen such that
$$
\int_{B(0,R_k)}\bigg[|\n u_k(x)|^2-\l\,\frac{|u_k(x)|^2}{|x|^2}\bigg]\,dx=
\int_{\R^N\setminus B(0,R_k)}\bigg[|\n u_k(x)|^2-\l\,\frac{|u_k(x)|^2}{|x|^2}\bigg]\,dx
=\frac12S_k(\l).
$$
Invariance by scaling yields
\begin{align}\label{eq:pid1}
\int_{\R^N}|\tilde u_k|^{2^*}=1,\quad &Q_{\l}(\tilde u_k)=S_k(\l),\\
\label{eq:pid2}\int_{\R^N}|v_k|^{2^*}=(S_k(\l))^{\frac N2},\quad &Q_{\l}(v_k)=(S_k(\l))^{\frac N2},
\end{align}
and
\begin{equation}\label{eq:sab1}
\int_{B(0,1)}\bigg[|\n v_k(x)|^2-\l\,\frac{|v_k(x)|^2}{|x|^2}\bigg]\,dx=
\int_{\R^N\setminus B(0,1)}\bigg[\n v_k(x)|^2-\l\,\frac{|v_k(x)|^2}{|x|^2}\bigg],dx
=\frac12(S_k(\l))^{\frac N2}.\hskip-.3cm
\end{equation}
Invariance by scaling also implies that $\{\tilde u_k\}_k$ is bounded in $\Di$, hence there exists a subsequence (still denoted as $\{\tilde u_k\}_k$) weakly converging to some $\tilde u_0$ in $\Di$. 

\medskip\noindent
{\bf Claim 1.} \quad We claim that $\tilde u_0\not\equiv 0$. 
Assume by contradiction that $\tilde u_0\equiv 0$. Using Proposition \ref{p:sab} for sequence $\tilde u_k$ with $r=1$ and $\e=\pm\frac14$ and taking into account \eqref{eq:sab1}, \eqref{eq:ven1}, and  \eqref{eq:ser1}, we deduce that there exist $\rho^+\in(0,1/4)$ and $\rho^-\in(-1/4,0)$ such that, up to a subsequence,
\begin{align}
\label{eq:ser2}&\int_{B(0,1+\rho^-)}|\n v_k|^2\to 0,\quad
\int_{B(0,1+\rho^-)}|v_k|^{2^*}\to 0,\quad\text{and }\int_{B(0,1+\rho^-)}\frac{|v_k|^2}{|x|^2}\to 0,\\
\label{eq:ser3}&
\int_{\R^N\setminus B(0,1+\rho^+)}|\n v_k|^2\to 0,\quad
\int_{\R^N\setminus B(0,1+\rho^+)}|v_k|^{2^*}\to 0,\quad\text{and }\int_{\R^N\setminus B(0,1+\rho^+)}\frac{|v_k|^2}{|x|^2}\to 0.
\end{align}
Note that weak convergence of $\tilde u_k\weakly 0$ in $\Di$, \eqref{eq:ven1}, and  \eqref{eq:ser1}, imply weak convergence of $v_k\weakly 0$ in $\Di$. Let $\eta$ be a smooth radial cut off function such that $0\leq\eta\leq1$, $\eta(x)\equiv 1$ for $1+\rho^-\leq|x|\leq 1+\rho^+$ and 
$\eta(x)\equiv 0$ for $|x|\not\in[3/4,5/4]$. Set $\tilde v_k:=\eta v_k$. Clearly $\tilde v_k\in\Dik$. By choice of $\eta$ and (\ref{eq:ser2}--\ref{eq:ser3}) we have
\begin{align*}
& Q_{\l}(\tilde v_k)=Q_{\l}(v_k)+o(1),\quad
\|\tilde v_k-v_k\|_{\Di}=o(1),\quad
\int_{\R^N}|\tilde v_k|^{2^*}=\int_{\R^N}|v_k|^{2^*}+o(1).
\end{align*}
Let us define 
$$
f(u):=\frac12\int_{B(0,5/4)\setminus B(0,3/4)}|\n u|^2-\frac1{2^*}
\int_{B(0,5/4)\setminus B(0,3/4)}|u|^{2^*},\quad u\in H^1_0(B(0,5/4)\setminus B(0,3/4)).
$$
From (\ref{eq:ser2}--\ref{eq:ser3}), \eqref{eq:ven1}, and \eqref{eq:pid2}, it is easy to verify that 
$$
f'(\tilde v_k)\to 0\quad\text{in }(\Di)^\star\quad\text{and}\quad f(\tilde v_k)=\frac1N (S_k(\l))^{N/2}+o(1)\leq \frac1N (S_{\rm circ}(\l))^{N/2}+o(1),
$$
i.e. $\tilde v_k$ is a Palais-Smale sequence for $f$ in
$H^1_0(B(0,5/4)\setminus B(0,3/4))$. From Struwe's representation
lemma for diverging Palais-Smale sequences \cite[Theorem
III.3.1]{Struwe}, we deduce the existence of an 
integer $M\in\N$, $M$ 
sequences of points $\{x^i_k\}_k\subset B(0,5/4)\setminus B(0,3/4)$
and $M$ sequences of radii $\{R^i_k\}_k$, $i=1,\dots,M$, such that
$\lim_kR^i_k=+\infty$ and  
\begin{equation}\label{eq:struwe}
\tilde v_k(x)=\sum_{i=1}^M (R^i_k)^{\frac{N-2}2} \tilde
v_0(R^i_k(x-x^i_k))+{\mathcal R}_k(x) \quad
 \text{where }{\mathcal 
  R}_k\to 0\text{ in }\Di, 
\end{equation}
$\tilde v_0=w_1^0$ and $w_1^0$ is the Talenti-Aubin function in
\eqref{eq:ZA}. 
Let us consider the sequence $\{x^1_k\}_k$; up to subsequence we can
assume that it converges to some point $x^1\subset\overline{
  B(0,5/4)\setminus B(0,3/4)}$. Let us write
$x^1=(z_1,y_1)\in\R^2\times\R^{N-2}$ where
$z_1=|z_1|e^{\bar\theta\sqrt{-1}}$. 

Let us first assume that $|z_1|\not=0$. We fix $J\in\N$ and for any $i=1,2,\dots,J$ we set
$$
S_i=\bigg\{(z,y)=(|z|e^{\theta\sqrt{-1}},y)\in \R^2\times\R^{N-2}:\ \bar\theta-\frac{(2i-1)\pi}J<\theta<\bar\theta+\frac{(2i+1)\pi}J\bigg\}.
$$
Note that $x^1\in S_1$ and there exists $\delta=\delta(J)>0$ such that 
$$B(x^1,\delta)\subset 
 \bigg\{(|z|e^{\theta\sqrt{-1}},y)\in \R^2\times\R^{N-2}:\ \bar\theta-\frac{(2i-1)\pi}{2J}<\theta<\bar\theta+\frac{(2i+1)\pi}{2J}\bigg\}.
$$
Choose $\bar k=\bar k(\delta)$ such that for all $k\geq\bar k$ 
$$
x^1_k=(z^1_k,y^1_k)\in B\Big(x^1,\frac\delta2\Big)\quad\text{and}\quad	 (R_k^1 )^{-1}<\frac\delta2.
$$
Moreover, if $\bar k$ is chosen sufficiently large, for each
$i=1,2,\dots,J$ it is possible to find $\tau_k^i\in \Z_k$ such that
$(\tau_k^i z_1,y_1)$ stays in the middle half of $S_i$, i.e.  
$$
(\tau_k^i z_1,y_1)\in 
\bigg\{(z,y)=(|z|e^{\theta\sqrt{-1}},y)\in 
\R^2\times\R^{N-2}:\ 
\bar\theta-\frac{(2i-1)\pi}{2J}<\theta<\bar\theta+
\frac{(2i+1)\pi}{2J}\bigg\}.
$$
Hence
$$ B\big((\tau_k^i z_1,y_1),\delta\big)\subset S_i,\quad (\tau_k^i z^1_k,y^1_k)\in 
B\Big((\tau_k^i z_1,y_1),\frac\delta2\Big)
$$
and consequently
$$ B\Big((\tau_k^i z_k^1,y_k^1),\frac\delta2\Big)\subset S_i
$$
which yields
$$ B\big((\tau_k^i z_k^1,y_k^1),(R_k^i)^{-1}\big)\subset S_i.
$$
In particular the $J$ balls $B\big((\tau_k^i z_k^1,y_k^1),(R_k^i)^{-1}\big)$ are disjoint, hence, by symmetry properties of $\tilde v_k$ we have that

\begin{align*}
(S_k(\l))^{N/2}+o(1)=\int_{\R^N}|\tilde v_k|^{2^*}&\geq \sum_{i=1}^J \int_{ B((\tau_k^i z_k^1,y_k^1),(R_k^i)^{-1})}|\tilde v_k|^{2^*}=\sum_{i=1}^J \int_{ B(x_k^1,(R_k^i)^{-1})}|\tilde v_k|^{2^*}.
\end{align*}
On the other hand from  \eqref{eq:struwe} we have, for $k$ large,
\begin{align*}
\bigg(&\int_{ B(x_k^1,(R_k^i)^{-1})}|\tilde v_k|^{2^*}\bigg)^{\frac1{2^*}}\\
&\geq 
\bigg(\int_{ B(x_k^1,(R_k^i)^{-1})}(R_k^1)^{N}\big|\tilde
v_0(R_k^1(x-x^1_k))\big|^{2^*}\bigg)^{\frac1{2^*}}- 
\bigg(\int_{ B(x_k^1,(R_k^i)^{-1})}|{\mathcal R}_k|^{2^*}\bigg)^{\frac1{2^*}}\geq\frac12 
\bigg(\int_{ B(0,1)}\tilde v_0^{2^*}\bigg)^{\frac1{2^*}}.
\end{align*}
Therefore
$$
(S_k(\l))^{N/2}+o(1)\geq \frac{J}{2^{2^*}}\int_{ B(0,1)}\tilde v_0^{2^*}
$$
and, in view of \eqref{eq:ven1}
$$
(S_{\rm circ}(\l))^{N/2}\geq \frac{J}{2^{2^*}}\int_{ B(0,1)}\tilde v_0^{2^*}.
$$
Letting $J\to+\infty$, we find a contradiction. Claim 1 is thereby proved in the case $|z_1|\not=0$. The case $|z_1|=0$ can be treated exploiting the radial symmetry of functions $\tilde u_k$ in the last $N-2$ variables with a similar argument (even simpler due to the stronger symmetry).

\medskip\noindent
{\bf Claim 2.} \quad We claim that $\tilde u_0\in\Dis$. 
We first note that $\tilde u_k$ satisfy the equation $-\Delta \tilde u_k-\l \frac{\tilde u_k}{|x|^2}=S_k(\l)\tilde u_k^{2^*-1}$. From \eqref{eq:ven1}, we can assume that $S_k(\l)\to L\in(0,+\infty)$ at least for a subsequence. Hence, due to weak convergence of $\tilde u_k\weakly \tilde u_0$,  we can pass to the limit in the equation to find that 
$\tilde u_0$ satisfies the equation $-\Delta \tilde u_0-\l \frac{\tilde u_0}{|x|^2}=L\tilde u_0^{2^*-1}$. By classical regularity theory for elliptic equations, we deduce that $\tilde u_0$ is a smooth function outside the origin.

Let $R>0$.  Assume that there exist $(z_1,y),(z_2,y)\in B(0,R)\cap(\R^2\times\R^{N-2})$, $|z_1|=|z_2|$, such that $\tilde u_0(z_1,y)\not=\tilde u_0(z_2,y)$. Then there exist $\delta>0$ such that $\tilde u_0(x)\not=\tilde u_0(y)$ for any $x\in B((z_1,y),\delta)$, $y\in B((z_2,y),\delta)$. 
Let $0<\e<\frac12 |B(0,\delta)|$. Since, up to a subsequence, $\tilde
u_k\to\tilde u_0$ a.e. in $B(0,R)$, by the Severini-Egorov Theorem,
there exists a measurable set $\Omega\subset B(0,R)$ such that $|\Omega|<\e$ and
$\tilde u_k\to\tilde u_0$ uniformly in $B(0,R)\setminus\Omega$. Hence
for $k$ large, $\tilde u_k(x)\not=\tilde u_k(y)$ for any $x\in
B((z_1,y),\delta)\setminus \Omega$, $y\in B((z_2,y),\delta)\setminus
\Omega$. On the other hand, if $k$ is large enough, there exists
$\tau_k\in\Z_k$ such that  
$$
\big|\tau_k(B((z_1,y),\delta))\triangle B((z_2,y),\delta)
\big|<\e,
$$
where $\triangle$ denotes the symmetric difference of sets.
Hence
$$
\big|\big(\tau_k(B((z_1,y),\delta))\cap B((z_2,y),\delta)\big)\setminus \Omega\big|>
|B(0,\delta)|-2\e>0.
$$
In particular the set $\big(\tau_k(B((z_1,y),\delta))\cap B((z_2,y),\delta)\big)\setminus \Omega$ has non-zero measure. If $(z,y)\in \big(\tau_k(B((z_1,y),\delta))\cap B((z_2,y),\delta)\big)\setminus \Omega$, then $z=\tau_k \tilde z$ with $(\tilde z,y)\in 
B((z_1,y),\delta)$ and by symmetry of $\tilde u_k$, $\tilde u_k(z,y)=\tilde u_k(\tilde z,y)$, thus giving a contradiction. Hence $\tilde u_0$ is invariant by the $\SO(2)$-action on the first two variables on $B(0,R)$ for any $R$. Invariance by the $\SO(N-2)$-action on the last $(N-2)$  variables follows easily from pointwise convergence.  
 Then we conclude that $\tilde u_0\in \Dis$.

\medskip\noindent
Hence we have proved that, up to a subsequence, $\tilde u_k\weakly\tilde u_0$ in $\Di$, with $\tilde u_0\in\Dis\setminus\{0\}$. Weak convegence yields
\begin{equation}\label{eq:pa1}
Q_{\l}(\tilde u_k)= Q_{\l}(\tilde u_0)+Q_{\l}(\tilde u_k-\tilde u_0)+o(1)
\end{equation}
while Brezis-Lieb Lemma implies
\begin{equation}\label{eq:pa2}
\|\tilde u_k\|_{2^*}^{2^*}= \|\tilde u_0\|_{2^*}^{2^*}+\|\tilde u_k-\tilde u_0\|_{2^*}^{2^*}+o(1).
\end{equation}
From \eqref{eq:pa1}, \eqref{eq:pa2}, \eqref{eq:pid1}, and \eqref{eq:minimizk} we have
$$
S_k(\l)\leq 
\frac{Q_{\l}(\tilde u_0)}{\|\tilde u_0\|_{2^*}^2} \leq S_k(\l)
\frac{Q_{\l}(\tilde u_k)-Q_{\l}(\tilde u_k-\tilde u_0)+o(1) }{\big(Q_{\l}(\tilde u_k)^{2^*/2}-Q_{\l}(\tilde u_k-\tilde u_0)^{2^*/2}+o(1)\big)^{2/2^*}}.
$$
Hence 
$$ 
\frac{Q_{\l}(\tilde u_k)-Q_{\l}(\tilde u_k-\tilde u_0)+o(1) }{\big(Q_{\l}(\tilde u_k)^{2^*/2}-Q_{\l}(\tilde u_k-\tilde u_0)^{2^*/2}+o(1)\big)^{2/2^*}}\geq 1.$$
Since $Q_{\l}(\tilde u_k)$ stay bounded away from $0$, this is possible only when $Q_{\l}(\tilde u_k-\tilde u_0)\to 0$. Since $Q_{\l}^{1/2}$ is an equivalent norm, we deduce that 
$\tilde u_k\to\tilde u_0$ in $\Di$. In particular $\|\tilde u_0\|_{2^*}=\lim_k\|\tilde u_k\|_{2^*}=1$.
Hence, by weakly lower semi-continuity of $Q_{\l}$, \eqref{eq:pid1}, and \eqref{eq:ven1}
$$
S_{\rm circ}(\l)\leq 
\frac{Q_{\l}(\tilde u_0)}{\|\tilde u_0\|_{2^*}^2} \leq \liminf_k  Q_{\l}(\tilde u_k)
=\liminf_k S_k(\l) \leq
\limsup_k  S_k(\l)\leq S_{\rm circ}(\l).
$$
Therefore all the above inequalities are indeed equalities. We have thus proved that along a subsequence, $S_k(\l)$ converges to $S_{\rm circ}(\l)$. The Uryson's property yields convergence of the entire sequence.~\end{pf}

\section{Proof of Theorem \ref{t:ach}}\label{sec:proof-theor-reft}

\noindent
The proof of Theorem \ref{t:ach} is based on Theorem \ref{t:ven} and the following lemma.

\begin{Lemma}\label{l:lis}
$\limsup_{k\to+\infty}S_k(\l_0,\L_1,\dots,\L_m)\leq S_{\rm circ}(\l_0,\L_1,\dots,\L_m)$.
\end{Lemma}
\begin{pf}
Let $\e>0$. Then from \eqref{eq:56} and density of ${\mathcal D}(\R^N\setminus\{0\})\cap\Dis$ in $\Dis$, there exists $u\in{\mathcal D}(\R^N\setminus\{0\})\cap \Dis$ such that $\int_{\R^N}|u|^{2^*}=1$ and
\begin{equation}\label{eq:pap1}
\int_{\R^N} |\n
      u|^2\,dx-{\l_0}\int_{\R^N}\frac{u^2(x)}{|x|^2}\,dx  -\sum_{\ell=1}^m
  {\L_{\ell}}\int_{\R^N}\bigg(\media_{S_{r_{\ell}}}\frac{u^2(y)}
{|x-y|^2}\,d\sigma(x)\bigg)\,dy< 
 S_{\rm circ}(\l_0,\L_1,\dots,\L_m)+\e.
\end{equation}
For any $\ell=1,\dots,m$, set 
$$
f_{\ell}(x):=\int_{\R^N}\frac{|u(y+x)|^2}{|y|^2}\,dy,\qquad x\in S_{r_{\ell}}.
$$
It is easy to check that $f_{\ell}\in C^0(S_{r_{\ell}})$; indeed if $x_n\in S_{r_{\ell}}$ converge to $x\in S_{r_{\ell}}$, by the Dominated Convergence Theorem we conclude that $\lim_n f_{\ell}(x_n)=f_{\ell}(x)$. Hence the Riemann sum
$$
\frac1k \sum_{i=1}^k f_{\ell}(a_i^{\ell})=\frac1k \sum_{i=1}^k
\int_{\R^N} \frac{|u(y)|^2}{|y-a_i^{\ell}|^2}\,dy
$$ 
converges to the integral
$$
\frac1{2\pi r_{\ell}}\int_{S_{r_{\ell}}}f_{\ell}(x)\,d\sigma(x)
=
\int_{\R^N}\bigg(\media_{S_{r_{\ell}}}\frac{u^2(y)}
{|x-y|^2}\,d\sigma(x)\bigg)\,dy.
$$
Hence there exists $\bar k$ such that for all $k\geq \bar k$ 
\begin{align}\label{eq:pap2}
\int_{\R^N} &|\n
      u|^2\,dx-{\l_0}\int_{\R^N}\frac{u^2(x)}{|x|^2}\,dx  -
\sum_{\ell=1}^m\sum_{i=1}^k
  \frac{\L_{\ell}}k\int_{\R^N}\frac{u^2(x)}
{|x-a_i^{\ell}|^2}\,dx-\e
\\
&\notag\leq
\int_{\R^N} |\n
      u|^2\,dx-{\l_0}\int_{\R^N}\frac{u^2(x)}{|x|^2}\,dx  -\sum_{\ell=1}^m
  {\L_{\ell}}\int_{\R^N}\bigg(\media_{S_{r_{\ell}}}\frac{u^2(y)}
{|x-y|^2}\,d\sigma(x)\bigg)\,dy.
\end{align}
From (\ref{eq:pap1}--\ref{eq:pap2}), we deduce that
$$
S_k(\l_0,\L_1,\dots,\L_m)-\e< S_{\rm circ}(\l_0,\L_1,\dots,\L_m)+\e.
$$
Taking $\limsup$ as $k\to+\infty$, since $\e$ is arbitrary we reach the conclusion.
\end{pf}

\begin{pfn}{Theorem \ref{t:ach}}
As in the proof of Theorem \ref{t:achcirc}, we can find a minimizing sequence $\{u_n\}_n\subset\Dik$ for
\eqref{eq:25} with the Palais-Smale property. Under assumption (\ref{eq:58}), Lemma \ref{l:comp-betw-conc} yields
\begin{equation}\label{eq:422}
S_{\rm circ}(\l_0,\L_1,\dots,\L_m)<
S_{\rm circ}(\l_0),
\end{equation}
while (\ref{eq:57}) implies 
\begin{equation}\label{eq:4222}
S_{k}(\l_0)\leq S_{k}\Big(\l_0+k\sum_{\ell=1}^k\l_{\ell}\Big).
\end{equation}
Let $0<\e<S_{\rm circ}(\l_0)-S_{\rm circ}(\l_0,\L_1,\dots,\L_m)$. From Lemma \ref{l:lis}, there exists $k_1=k_1(\e)$ such that for all $k\geq k_1$
\begin{equation}\label{eq:423}
S_{k}(\l_0,\L_1,\dots,\L_m)<
S_{\rm circ}(\l_0,\L_1,\dots,\L_m)+\e.
\end{equation}
From Theorem \ref{t:ven} \and \eqref{eq:422},  there exists $k_2$ such that for all $k\geq k_2$
\begin{equation}\label{eq:424}
S_{\rm circ}(\l_0,\L_1,\dots,\L_m)+\e<S_{k}(\l_0)\leq S_{\rm circ}(\l_0).
\end{equation}
Let $k_3$ be such that for all $k\geq k_3$
\begin{equation}\label{eq:425}
S_{\rm circ}(\l_0,\L_1,\dots,\L_m)+\e\leq
\min\bigg\{k^{\frac2N}S,
k^{\frac2N}S(\l_1), \dots, k^{\frac2N}S(\l_m)\bigg\}.
\end{equation}
From (\ref{eq:4222}--\ref{eq:425}), we conclude that 
 for all $k\geq \max\{k_1,k_2,k_3\}$ 
$$
S_k(\l_0,\L_1,\dots,\L_m)<\min\bigg\{k^{\frac2N}S,
k^{\frac2N}S(\l_1), \dots, k^{\frac2N}S(\l_m), S_k(\l_0),
S_k\Big(\l_0+k\sum\nolimits_{\ell=1}^{m}\l_{\ell}\Big)\bigg\}.
$$
From above and the Palais-Smale condition proved in
Theorem \ref{th:ps}, we deduce that $\{u_n\}_{n\in\N}$ has a  subsequence strongly
converging to some $u_0\in\Dik$ such that 
$J_k(u_0)=\frac1NS_k(\l_0,\L_1,\dots,\L_m)$.
Hence $u_0$ achieves the infimum in \eqref{eq:25}. Since
$J_k$ is even, also $|u_0|$ is a minimizer in \eqref{eq:25}
and then $v_0=S_k(\l_0,\L_1,\dots,\L_m)^{1/(2^*-2)}|u_0|$ is a
nonnegative solution to equation \eqref{eq:22}. The maximum principle
implies the positivity  outside singularities of such a solution.~\end{pfn}

\section{Proof of Theorem \ref{t:ach1}}\label{sec:proof-theor-reft-1}

\noindent We now provide a sufficient condition for the
infimum 
in (\ref{eq:25}) to stay below the level $k^{2/N}S(\l_j)$, in
correspondence of   which 
possible concentration at singular points located at the $j$-th
polygon can occur.
We denote by $\Theta_{j\ell}$ the minimum angle formed by vectors
$a_i^j$ and $a_s^{\ell}$, see figure below.
\begin{center}
 \vskip0.5truecm\noindent
 \epsfxsize=1.5in \epsfbox{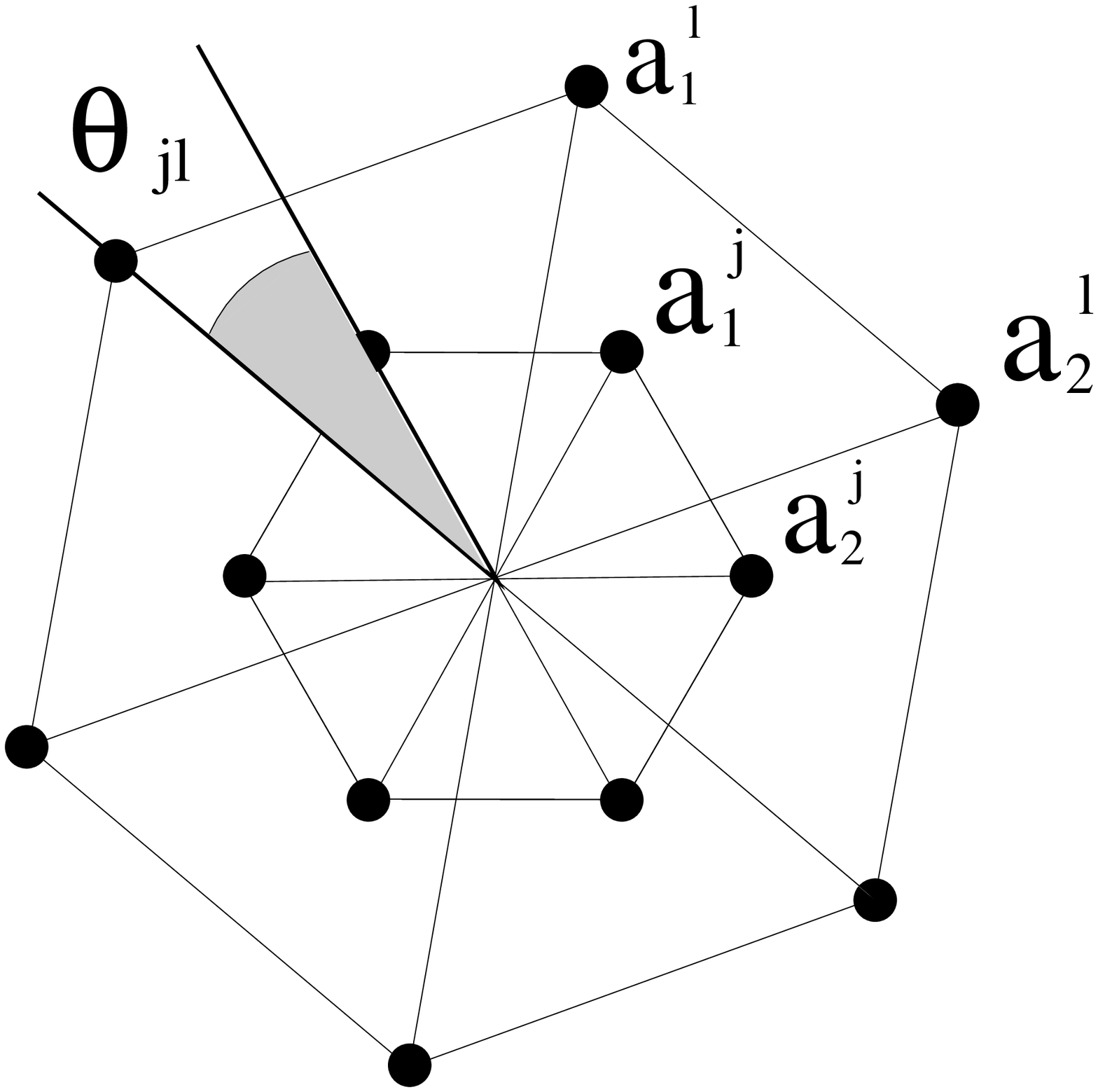}\\
{\scriptsize Figure 4 (The angle $\Theta_{j\ell}$.)}
\end{center}
The following lemma can be proved by standard trigonometry calculus.
\begin{Lemma}\label{l:trigo}
For any $i,s=1,2,\dots,k$, and $j,\ell=1,2,\dots,m$, there holds
\begin{align*}
& |a_i^{j}-a_s^j|=2r_j\Big|\sin\frac{(s-i)\pi}k\Big|,\\
& | a_i^j-a_s^{\ell}|^2=r_j^2+r_{\ell}^2-2 r_j
r_{\ell}\cos\Big(\frac{2\pi(i-s)}k+\Theta_{j\ell}\Big).
\end{align*} 
\end{Lemma}

 \begin{Lemma}\label{cor:sub}
Let $j\in\{1,2,\dots,m\}$. If 
\begin{align}
\label{eq:37}
0<\l_j\leq\frac{N(N-4)}4
\end{align}
and
\begin{align}\label{eq:50}
\frac{\l_0}{|r_j|^2}+\l_j
\sum_{i=1}^{k-1}\frac{1}
{4r_j^2\big|\sin\frac{i\pi}k\big|^2}+\sum_{\substack{\ell=1\\\ell\not
    =j}}^m\l_{\ell}\sum_{i=1}^k\frac{1}{r_j^2+r_{\ell}^2-2 r_j
r_{\ell}\cos\Big(\frac{2\pi i}k+\Theta_{j\ell}\Big)
}>0,
\end{align}
then
\begin{align}\label{eq:16}
S_k(\l_0,\L_1,\dots,\L_m)<k^{2/N}S(\l_j).
\end{align}
\end{Lemma}
\begin{pf}
Let $z(x)=\sum_{i=1}^k z_{\mu}^{\l_j}(x-a_i^j)\in \Dik$. Then
\begin{align*}
S_k(\l_0,\L_1,\dots,\L_m)&\leq 
\frac{Q_{\l_0,\l_1,\dots,\l_m}(z)}{{\displaystyle{ 
    \bigg(\int_{\R^N}  |z|^{2^*}dx\bigg)^{2/2^*}}}},
\end{align*}
where $Q_{\l_0,\l_1,\dots,\l_m}$ denotes the quadratic form defined by
$$
Q_{\l_0,\l_1,\dots,\l_m}(u)
={\displaystyle{\int_{\R^N} |\n
      u|^2dx-\int_{\R^N}\frac{\l_0}{|x|^2}u^2(x)\,dx 
-\sum_{\ell=1}^m\sum_{i=1}^k 
  {\l_{\ell}}\int_{\R^N}\frac{u^2(x)}
{|x-a_i^{\ell}|^2}\,dx}}.
$$
Note that 
$$
 \bigg(\int_{\R^N}  |z|^{2^*}dx\bigg)^{2/2^*}\geq k^{2/2^*}.
$$
Moreover, from (\ref{eq:minimiz}), Lemmas \ref{l:interest} and \ref{l:gradgrad}
we find that
\begin{align*}
&Q_{\l_0,\l_1,\dots,\l_m}(z)=k\int_{\R^N}|\n z_{\mu}^{\l_j}|^2\,dx-k \l_j
\int_{\R^N}\frac{|z_1^{\l_j}|^2}{|x|^2}\,dx-\l_0\sum_{i=1}^k\int_{\R^N}
\frac{|z_{\mu}^{\l_j}|^2}{|x+a_i^j|^2}\,dx\\
&-\sum_{\ell=1}^m\sum_{\substack{i=1,s=1\\(i,\ell)\not=(s,j)}}^k\int_{\R^N}
\frac{\l_{\ell}|z_{\mu}^{\l_j}|^2}{|x+a_s^j-a_i^{\ell}|^2}\,dx
+\sum_{\substack{i=1,s=1\\i\not=s}}^k\int_{\R^N}\n
z_{\mu}^{\l_j}(x-a_i^j)\cdot\n z_{\mu}^{\l_j}(x-a_s^j)\,dx\\
&-\sum_{\ell=1}^m\sum_{\substack{i=1,s=1, t=1\\s\not=t}}^k\int_{\R^N}
\frac{\l_{\ell}z_{\mu}^{\l_j}(x-a_s^j)z_{\mu}^{\l_j}(x-a_t^j)}
{|x-a_i^{\ell}|^2}\,dx
-\l_0\sum_{\substack{i=1,s=1\\i\not=s}}^k
\int_{\R^N}
\frac{z_{\mu}^{\l_j}(x-a_i^j)z_{\mu}^{\l_j}(x-a_s^j)}
{|x|^2}\,dx\\
&=
\begin{cases}
k S(\l_j)-\mu^2\Big(\int_{\R^N}|z_1^{\l_j}|^2\Big)\bigg[
\frac{k\l_0}{|r_j|^2}+{\displaystyle{\sum_{\ell=1}^m\sum_{\substack{i=1,s=1
\\(i,\ell)\not=(s,j)}}^k}}\frac{\l_{\ell}}{|a_i^{\ell}-a_s^j|^2}+o(1)\bigg]
\text{if }\l_j<\frac{N(N-4)}4, \\
k S(\l_j)-\mu^2|\ln \mu|\a_{\l_j,N}^2\bigg[
\frac{k\l_0}{|r_j|^2}+{\displaystyle{\sum_{\ell=1}^m\sum_{\substack{i=1,s=1
\\(i,\ell)\not=(s,j)}}^k}}\frac{\l_{\ell}}{|a_i^{\ell}-a_s^j|^2}+o(1)\bigg]
\text{if }\l_j=\frac{N(N-4)}4.
\end{cases}
\end{align*}
Therefore (\ref{eq:16}) holds provided
\begin{align}\label{eq:36}
\frac{k\l_0}{|r_j|^2}+{\displaystyle{\sum_{\ell=1}^m\sum_{\substack{i=1,s=1
\\(i,\ell)\not=(s,j)}}^k}}\frac{\l_{\ell}}{|a_i^{\ell}-a_s^j|^2}>0 .
\end{align}
It is easy to verify that assumption (\ref{eq:50}) and Lemma
\ref{l:trigo}  imply (\ref{eq:36}).
\end{pf}

\begin{Lemma}\label{l:largek}
For any $j,\ell=1,\dots,m$, $j\not=\ell$,  there holds
\begin{align}\label{eq:51}
&\lim_{k\to\infty}\frac1{\Lambda_j}\bigg[\l_j
\sum_{i=1}^{k-1}\frac{1}
{4r_j^2\big|\sin\frac{i\pi}k\big|^2}\bigg]=+\infty,\\
&\label{eq:52}
\l_{\ell}\sum_{i=1}^k\frac{1}{r_j^2+r_{\ell}^2-2 r_j
r_{\ell}\cos\Big(\frac{2\pi
  i}k+\Theta_{j\ell}\Big)}=\Lambda_{\ell}\, O(1)\quad\text{as
}k\to+\infty.
\end{align}
\end{Lemma}
\begin{pf}
A direct calculation yields
\begin{align*}
\frac{\l_j}{\Lambda_j}
\sum_{i=1}^{k-1}\frac{1}
{4r_j^2\big|\sin\frac{i\pi}k\big|^2}&\geq 
\frac 1k \int_1^{k/2}\frac{ds}{4r_j^2\big|\sin\frac{s\pi}k\big|^2}\geq
\frac 1\pi \int_{\pi/k}^{\pi/2}\frac{dt}{4r_j^2|\sin t|^2}
\mathop{\longrightarrow}\limits_{k\to+\infty}+\infty
.
\end{align*}
On the other hand
\begin{align*}
\frac{\l_{\ell}}{\Lambda_{\ell}}\sum_{i=1}^k\frac{1}{r_j^2+r_{\ell}^2-2 r_j
r_{\ell}\cos\Big(\frac{2\pi
  i}k+\Theta_{j\ell}\Big)}&\leq \frac1k\int_0^{k} 
\frac{ds}{r_j^2+r_{\ell}^2-2 r_j
r_{\ell}\cos\Big(\frac{2\pi
  s}k+\Theta_{j\ell}\Big)}\\
&= \frac1{2\pi}\int_0^{2\pi} 
\frac{dt}{r_j^2+r_{\ell}^2-2 r_j
r_{\ell}\cos(t+\Theta_{j\ell})}\in\R^+,
\end{align*}
thus proving (\ref{eq:52}).
\end{pf}

\begin{remark}\label{rem:klrage}
Lemma \ref{l:largek} implies that if we fix
$\l_0,\Lambda_1,\dots,\Lambda_m$ and let $k\to+\infty$, then the
quantity in formula (\ref{eq:50}) tends to $+\infty$. Hence condition 
(\ref{eq:50}) is satisfied for $k$ sufficiently large.
\end{remark}

\noindent Let us now compare levels $S_k(\l_0,\L_1,\dots,\L_m)$ and
$S_k(\l_0)$ which is related to concentration at the origin. Two cases
can occur:
\begin{align*}
{\it (i)}\hskip1cm&S_k(\l_0)\geq k^{2/N}S,\\
{\it (ii)}\hskip1cm&S_k(\l_0)< k^{2/N}S.
\end{align*}
In case  (i), since $S=S(0)$ and $\lambda\mapsto S(\l)$ is a
nonincreasing function, to exclude that the infimum in
(\ref{eq:minimizk}) stays above $S_k(\l_0)$ it is enough to compare
$S_k(\l_0,\L_1,\dots,\L_m)$ with  $k^{2/N}S(\l_j)$ where
$\l_j=\max\{\l_{\ell}\}_{1\leq\ell\leq m}$, as we have done in Lemma
\ref{cor:sub}.

The study of  case (ii) is based on Lemma \ref{l:ska}. Indeed, using  
Lemma \ref{l:ska} and estimates of Lemma
\ref{l:interest}, we can 
prove the following lemma. 

\begin{Lemma}\label{l:sub}
If $N\geq 4$, $S_k(\l_0)<k^{2/N}S$, and one of the following
assumptions is satisfied 
\begin{align}
\label{eq:38}
\l_0\leq\frac{N(N-4)}4\quad&\text{and}\quad
\sum_{\ell=1}^m\frac{\l_{\ell}}{|r_{\ell}|^2}>0,\\
\label{eq:39}
\frac{N(N-4)}4<\l_0<\frac{(N-2)^2}4\quad&\text{and}\quad\sum_{\ell=1}^m
\frac{\l_{\ell}}
{|r_{\ell}|^{\sqrt{(N-2)^2-4\l_0}}}>0,
\end{align}
then
\begin{equation}\label{eq:41}
S_k(\l_0,\L_1,\dots,\L_m)<S_k(\l_0).
\end{equation}
\end{Lemma}
\begin{pf}
From Lemma \ref{l:ska}, we have that $S_k(\l_0)$ is attained by some
$u^{\l_0}\in\Dik$. By homogeneity of the Rayleigh quotient, we can
assume $\int |u^{\l_0}|^{2^*}=1$. Furthermore, the function
  $v^{\l_0}=S_k(\l_0,\L_1,\dots,\L_m)^{1/(2^*-2)}|u^{\l_0}|$ is a
  nonnegative solution to (\ref{eq:1}), hence we can apply
  Lemma~\ref{l:interest} to study the behavior of 
  $\int_{\R^N}\frac{|u^{\l_0}_{\mu}|^2}{|x+\xi|^2}\,dx$ as $\mu\to
  0$. Hence we obtain that there exists some positive constant
  $\kappa_0$ such that
\begin{align}\label{eq:71}
&S_k(\l_0,\L_1,\dots,\L_m)\\
&\notag \leq \int_{\R^N}|\n
u_{\mu}^{\l_0}(x)|^2\,dx-\l_0\int_{\R^N}\frac{|
  u_{\mu}^{\l_0}(x)|^2}{|x|^2}\,dx -
\sum_{\ell=1}^m\sum_{i=1}^k\l_{\ell} \int_{\R^N}\frac{|
  u_{\mu}^{\l_0}(x)|^2}{|x-a_i^{\ell}|^2}\,dx\\
\notag &=S_k(\l_0)-
\begin{cases}
\mu^2\int_{\R^N}|u^{\l_0}|^2\left({\displaystyle{\sum_{\ell=1}^m\sum_{i=1}^k}}\dfrac{\l_{\ell}}{|r_{\ell}|^2}+o(1)\right)
    &\text{if }\l_0<\frac{N(N-4)}4\\[15pt]
\kappa_0^2\mu^2|\ln\mu|\left({\displaystyle{\sum_{\ell=1}^m
\sum_{i=1}^k}}\dfrac{\l_{\ell}}{|r_{\ell}|^2}
+o(1)\right)
    &\text{if }\l_0=\frac{N(N-4)}4\\[15pt]
\kappa_{0}^2\,\beta_{\l,N}\,
\mu^{\sqrt{(N-2)^2-4\lambda}}
\left({\displaystyle{\sum_{\ell=1}^m
\sum_{i=1}^k}}\dfrac{\l_{\ell}}{|r_{\ell}|^{\sqrt{(N-2)^2-4\lambda_0}}}
+o(1)\right)
    &\text{if }\l_0>\frac{N(N-4)}4.
\end{cases}
\end{align}
Taking $\mu$ sufficiently small we obtain that
either assumption (\ref{eq:38}) or (\ref{eq:39}) yield
(\ref{eq:41}).
\end{pf}

\begin{pfn}{Theorem \ref{t:ach1}}
As in the proof of Theorem \ref{t:achcirc}, we can find a minimizing
sequence $\{u_n\}_{n\in\N}$ which has the Palais-Smale property, more precisely  $J_k'(u_n)\to 0$ in $(\Dik)^{\star}$ and 
$J_k(u_n)\to \frac1N S_k(\l_0\L_1,\dots, \L_m)$. Assumption (\ref{eq:42}) yields
\begin{equation}\label{eq:48}
S_k\Big(\l_0+k\sum_{\ell=1}^k\l_{\ell}\Big)\geq S_k(\l_0).
\end{equation}
Note also that (\ref{eq:43}) and (\ref{eq:47}) imply that $N> 4$
and $\l_m>0$.   
Two cases can occur. If $S_k(\l_0)<k^{\frac2N}S$, then Lemma
\ref{l:sub} and assumption (\ref{eq:47})
yield
\begin{equation*}
 S_k(\l_0,\L_1,\dots,\L_m)<S_k(\l_0).
\end{equation*}
If $S_k(\l_0)\geq  k^{2/N}S$, from monotonicity we have $S_k(\l_0)\geq
k^{2/N}S>k^{2/N}S(\l_m)$. In both cases from Lemma \ref{cor:sub} and
(\ref{eq:46}) we deduce
$k^{2/N}S(\l_m)>S_k(\l_0,\L_1,\dots,\L_m)$. Hence we obtain that 
$$
S_k(\l_0,\L_1,\dots,\L_m)<\min\bigg\{k^{\frac2N}S,
k^{\frac2N}S(\l_1), \dots, k^{\frac2N}S(\l_m), S_k(\l_0),
S_k\Big(\l_0+k\sum\nolimits_{\ell=1}^{m}\l_{\ell}\Big)\bigg\}.
$$
From above and the Palais-Smale condition proved in
Theorem \ref{th:ps}, we deduce that $\{u_n\}_{n\in\N}$ has a  subsequence strongly
converging to some $u_0\in\Dik$ 
which achieves the infimum in
\eqref{eq:25}. Moreover $v_0=S_{k}(\l_0\L_1,\dots,\L_m)^{\frac1{2^*-2}}|u_0|$ is a solution to \eqref{eq:22}.
\end{pfn}

\begin{remark}
Theorem  \ref{t:ach1} contains an alternative proof to Theorem \ref{t:ach} in the case $N>4$, as it follows easily gathering Theorem  \ref{t:ach1} and Remark
\ref{rem:klrage}. Note that the assumption $N>4$ is needed to ensure
that (\ref{eq:43}) and (\ref{eq:47}) hold. 
However, with respect to Theorem \ref{t:ach}, it contains a more precise information on how $k$ must be large in order to solve the problem.
\end{remark}

\section{Proof of Theorem
  \ref{t:posmas}}\label{sec:proof-theor-reft:p}

\begin{pfn}{Theorem \ref{t:posmas}} Assume first (i).
Let $\e>0$.
Then from \eqref{eq:25} and density of ${\mathcal
  D}(\R^N\setminus\{0\})\cap\Dik$ in $\Dik$, there exists
$u\in{\mathcal D}(\R^N\setminus\{0\})\cap \Dik$ such that
$Q_{\l_0}(u)\leq S_k(\l_0)+\e$. Let 
$$
u_{\mu}(x)=\mu^{-\frac{N-2}2}u(x/\mu), \quad\mu>0.
$$
By Dominated Convergence Theorem it is easy to verify that 
$$
\lim_{\mu\to0}\int_{\R^N}\frac{|u_{\mu}(x)|^2}{|x-a_i^{\ell}|^2}=0,
$$
hence
\begin{align*}
S_k(\l_0,\L_1,\dots, \L_m)&\leq
\frac{{{\int_{\R^N} |\n
      u_{\mu}|^2dx-\l_0\int_{\R^N}
\frac{u_{\mu}^2(x)}{|x|^2}-
\int_{\R^N}\sum_{\ell=1}^m\frac{\L_{\ell}}k\sum_{i=1}^k
  \frac{u_{\mu}^2(x)}{|x-a_i^{\ell}|^2}
dx}}}
{{{ 
    \big(\int_{\R^N}  |u_{\mu}|^{2^*}dx\big)^{2/2^*}}}}\\[10pt]
&=Q_{\l_0}(u)+o(1)\leq  S_k(\l_0)+\e+o(1)\quad\text{as }\mu\to0.
\end{align*}
Letting $\mu\to0$ we have that $S_k(\l_0,\L_1,\dots, \L_m)\leq
S_k(\l_0)+\e$ for all $\e>0$. Hence
$$
S_k(\l_0,\L_1,\dots, \L_m)\leq
S_k(\l_0).
$$
 Assume by
contradiction that $ S_{k}(\l_0,\L_1,\dots,\L_m)$ is attained by  $\bar u\in
\Dik\setminus\{0\}$, then
\begin{align*}
S_{k}(\l_0,\L_1,\dots,
\L_m)&=
\frac{{{\int_{\R^N} |\n
      \bar u|^2dx-\l_0\int_{\R^N}
\frac{\bar u^2(x)}{|x|^2}-
\int_{\R^N}\sum_{\ell=1}^m\frac{\L_{\ell}}k\sum_{i=1}^k
  \frac{\bar u^2(x)}{|x-a_i^{\ell}|^2}
dx}}}
{{{ 
    \big(\int_{\R^N}  |\bar u|^{2^*}dx\big)^{2/2^*}}}}
\\
&>
\frac{{ {\int_{\R^N} |\n
     \bar  u|^2\,dx-\int_{\R^N}\frac{\l_0}{|y|^2}\bar u^2(y)\,dy
}}}{{ { 
    \big(\int_{\R^N}  |\bar u|^{2^*}dx\big)^{2/2^*}}}}
\geq S_{k}(\l_0),
\end{align*}
giving rise to a contradiction. The proof of non-attainability of 
$S_{\rm circ}(\l_0,\L_1,\dots,\L_m)$ is analogous and is based on
\eqref{eq:70}.  

Assume now that (ii) holds. Then for all $u\in\Dis$, $u\geq 0$,
denoting by   $u^*$  the Schwarz symmetrization of $u$ (see
\eqref{eq:82}), from \eqref{eq:wil} and \eqref{eq:wil1} it follows
that
\begin{align*}
\frac{Q^{\rm circ}_{\l_0,\L_1,\dots,\L_m}(u)}{\big(\int_{\R^N}  |
  u|^{2^*}dx\big)^{2/2^*}}&\geq 
\frac{\int_{\R^N} |\n
      u^*|^2\,dx-\big(\l_0+
\sum_{\ell=1}^m
  {\L_{\ell}}\big)\int_{\R^N}\frac{|u^*(y)|^2}{|y|^2}\,dy}
{\big(\int_{\R^N}  |
  u^*|^{2^*}dx\big)^{2/2^*}}
\\
&\geq S_{\rm circ}\Big(\l_0+\sum_{\ell=1}^m {\L_{\ell}}\Big)=
S\Big(\l_0+\sum_{\ell=1}^m {\L_{\ell}}\Big).
\end{align*}
Hence $S_{\rm circ}(\l_0,\L_1,\dots,\L_m)\geq  S_{\rm
  circ}\big(\l_0+\sum_{\ell=1}^m {\L_{\ell}}\big)$. On the other hand,
setting $\Lambda=\l_0+\sum_{\ell=1}^m {\L_{\ell}}$, and using
\cite[Corollary 3.2]{FT} we obtain
\begin{align*}
S_{\rm circ}(\l_0,\L_1,\dots,\L_m)&\leq 
Q^{\rm circ}_{\l_0,\L_1,\dots,\L_m}(z_{\mu}^{\Lambda})=
Q_{\Lambda}(z_{\mu}^{\Lambda})+o(1)=S(\Lambda)+o(1)\\
&=S_{\rm circ}\Big(\l_0+\sum_{\ell=1}^m
{\L_{\ell}}\Big)+o(1)\quad\text{as }\mu\to\infty.
\end{align*}
Then 
\begin{equation}\label{eq:83}
S_{\rm circ}(\l_0,\L_1,\dots,\L_m)=  S_{\rm
  circ}\Big(\l_0+\sum_{\ell=1}^m {\L_{\ell}}\Big).
\end{equation}
If $S_{\rm circ}(\l_0,\L_1,\dots,\L_m)$ was attained by some $\bar
u\in\Dis\setminus\{0\}$, then 
$$
S_{\rm circ}(\l_0,\L_1,\dots,\L_m)=\frac{Q^{\rm
    circ}_{\l_0,\L_1,\dots,\L_m}(\bar u)}{\big(\int_{\R^N}  |
 \bar u|^{2^*}dx\big)^{2/2^*}}\geq 
\frac{Q_{\Lambda}(\bar u^*)}
{\big(\int_{\R^N}  |
  \bar u^*|^{2^*}dx\big)^{2/2^*}}
\geq S_{\rm circ}(\Lambda).
$$
Due to (\ref{eq:83}), all above inequalities are indeed equalities; in
particular $Q^{\rm
    circ}_{\l_0,\L_1,\dots,\L_m}(\bar u)=Q_{\Lambda}(\bar u^*)$ which
    (taking into account  Polya-Szego inequality and \eqref{eq:wil})
    yields
\begin{align*}
0&\leq \int_{\R^N} |\n \bar u|^2dx-\int_{\R^N} |\n  \bar u^*|^2dx\\&=
\int_{\R^N}\bigg(\frac{\l_0}{|y|^2}+
\sum_{\ell=1}^m
  {\L_{\ell}}\media_{S_{r_{\ell}}}\frac{d\sigma(x)}
{|x-y|^2}\bigg)\bar u^2(y)\,dy-\Lambda\int_{\R^N}\frac{|\bar
  u^*(y)|^2}{|y|^2}\,dy\leq 0.
\end{align*}
Then $\int_{\R^N} |\n \bar u|^2dx=\int_{\R^N} |\n  \bar
u^*|^2dx$. From \cite{BZ}, it follows that  $\bar u$ must be spherically symmetric with
respect to  some point. Since $\bar u$  is a solution
to equation \eqref{eq:55} (up to some Lagrange's multiplier), the
potential in equation \eqref{eq:55} must be spherically symmetric,
thus giving rise to a contradiction. The proof of non-attainability of 
$S_{k}(\l_0,\L_1,\dots,\L_m)$ is contained in \cite[Theorem 1.3]{FT}.
\end{pfn}



\end{document}